%% file: sphere_LB.tex
\documentclass[11pt]{article}
\usepackage{latexsym,amsmath, mathrsfs,amssymb,tabu,
	CJKutf8,graphicx,bm,blkarray,delarray,float,color,soul,extarrows, ulem, framed}
	\usepackage[thicklines]{cancel}
\newtheorem{proof}{\it Proof}

\usepackage{epsfig,graphics,subfigure,psfrag,amsmath,amssymb}
\usepackage{enumerate}
\usepackage{bm}
\usepackage{graphicx}
\usepackage{epsfig}
\usepackage{epstopdf}
\usepackage{subfigure}
\usepackage{caption}
\usepackage{overpic}
\usepackage{epsfig}
\usepackage{xcolor}
\usepackage{sidecap}
\usepackage{fancyhdr}
\usepackage{color,tabu,diagbox,multirow,caption}
\usepackage{algorithm} 
\usepackage{algorithmic} 
\usepackage{float}
\usepackage{fullpage} 
\usepackage{xkeyval}
\usepackage[draft]{changes}
\usepackage{cite}

\newcommand{\bbS}{\mathbb{S}}

\newcommand{\bbI}{\mathbb{I}}
\newcommand{\bbZ}{\mathbb{Z}}

\newcommand{\mbk}{\bm{k}}

\renewcommand{\Re}{\operatorname{Re}}
\renewcommand{\Im}{\operatorname{Im}}

\newcommand{\FR}[1]{{\color{black}{#1}}}
\newcommand{\TR}[1]{{\color{black}{#1}}}
\newcommand{\tR}[1]{{\color{black}{#1}}}
\newcommand{\fR}[1]{{\color{black}{#1}}}
\usepackage[utf8]{inputenc}
\newcounter{Algctr}[section]
\renewcommand{\theAlgctr}{\thesection.\arabic{Algctr}}

\newtheorem{theorem}{Theorem}[section]

\newtheorem{remark}{Remark}[section]

\numberwithin{equation}{section}

\begin{document}
\setulcolor{red} 
\setstcolor{red} 
\sethlcolor{yellow} 

\title{\bf Efficient numerical methods for computing stationary states of spherical Landau-Brazovskii model}
\author{Qun Qiu$^{1}$, Wei Si$^{2}$,Guanghua Ji$^{1}$\;$^{*}$ and Kai Jiang$^{3}$
\footnote{Corresponding author. \\Emails: ghji@bnu.edu.cn, kaijiang@xtu.edu.cn} }
     \date{}
 \maketitle
\begin{center}
\small $^{1}$School of Mathematical Sciences, Beijing Normal University,Beijing, China, 100875\\
\small $^{2}$ School of Mathematics Sciences, Peking University, Beijing, China, 100871 \\
\small $^{3}$School of Mathematics and Computational Science, Hunan Key Laboratory for Computation and Simulation in Science and Engineering, Key Laboratory of Intelligent Computing and  Information Processing of Ministry of Education,  Xiangtan University, Xiangtan, Hunan, China, 411105 \\
\end{center}
\input{Abstract}

\input{Part_1_Introduction}

\input{Part_2_Problem_formulation}

\input{Part_3_Numerical_methods}

\input{Part_4_Numerical_results}

\input{Conclusion}

\nocite{*}
\normalem


\end{document}

%% file: Abstract.tex
\begin{abstract}
In this paper, we develop a set of efficient methods to compute stationary states of the spherical Landau-Brazovskii (LB) model in a discretization-then-optimization way.
		First, we discretize the spherical LB energy functional into a finite-dimensional energy function by the spherical harmonic expansion.
		Then five optimization methods are developed to compute stationary states of the discretized energy function, including the accelerated adaptive Bregman proximal gradient, Nesterov, adaptive Nesterov, adaptive nonlinear conjugate gradient and adaptive gradient descent methods.
		To speed up the convergence, we propose a principal mode analysis (PMA) method to estimate good initial configurations and sphere radius.
		The PMA method also reveals the relationship between the optimal sphere radius and the dominant degree of spherical harmonics. 
		Numerical experiments show that our approaches significantly reduce the number of iterations and the computational time.
\end{abstract}

%% file: Part_1_Introduction.tex
	\section{Introduction}
	\label{sec1}
	\label{sec:Intro}
	Landau models are powerful tools for studying the microscopic behavior of structures in physics and materials science, such as symmetry breaking
	\cite{1950_LG}, pattern formation \cite{1975Brazovskii,1977_SH,Zhang2008}
	and phase transitions \cite{Shi2003,Zhanghao2014}.
	These models utilize order parameter functions to characterize the degree of order in the system.
	One specific Landau model of interest is the Landau-Brazovskii (LB) model \cite{1975Brazovskii}, which has proven valuable in describing periodic crystals and phase transitions in Euclidean space \cite{Zhang2008,Shi2003,Zhanghao2014}.
	Recently, the spherical LB model has been widely employed to explore pattern formation \cite{pattern_forming_2003,spiral_symmetry_2010}, block copolymer assembly \cite{Zhang_2014}, and viral capsids
	\cite{SLB_Sanjay_2017} on a spherical surface.
	Compared with the Swift-Hohenberg model without three-body interaction terms \cite{1977_SH,2002_Elder_SH, 2010_Backofen_Particles, 2011_Backofen_Continuous, 2012_Aland_Particles,2011_Aland_Continuum,2012_Aland_Bucking, 2018_Praetorius_Active}, the LB model can describe the first-order phase transition \cite{1975Brazovskii,1987Fredrickson,Zhanghao2014}.
	In this work, we focus on the development of efficient methods of finding equilibrium ordered structures of spherical LB model instead of studying quasi-equilibrium dynamical phase behavior.
	
	\par
	The presence of multiple solutions and non-linearity poses a \FR{challenge} in designing fast and efficient methods for quickly finding stationary states.
	Efficient computation of stationary states of the Landau free energy functional, corresponding to ordered structures, is essential due to their significance in determining material properties.
	Generally, existing numerical approaches for computing stationary states of the
	Landau free energy functional can be divided into three categories.
	The first category involves solving the Euler-Lagrange equation, the first-order variation of free energy functional.
	The second category comprises gradient flow approaches, such as general semi-implicit methods \cite{2012_Aland_Bucking,Zhang_2014,2018_Praetorius_Active}, stabilized factor 
	methods \cite{SIS_Shen_2010,CN_Feng_2013}, exponential time difference 
	schemes \cite{spiral_symmetry_2010}, convex splitting \cite{CSplitting_Vignal_2015,CSplitting_Lee_2016},
	operator splitting \cite{OSplitting_Lee_2015,OSplitting_Zhai_2021} and
	auxiliary variable methods \cite{IEQ_Yang_2016,CN_Yang_2017,SAV_Shen_2019}.
	Except the time discretization methods, spatial discretization techniques are important to numerically \FR{solve} gradient flow equations, including the finite difference method \cite{2022_Yang_FDM}, the finite volume method \cite{2005_Tang_Phase}, the finite element method \cite{2010_Backofen_Particles,2012_Aland_Bucking,2012_Aland_Particles} and the spectral method \cite{Zhang_2014,2018_Praetorius_Active}.
	Gradient flow approaches primarily focus on the dynamic evolution of ordered structures.
	The third approach treats the problem as an optimization task, directly computing stationary states of free energy
	functional using optimization algorithms based on proper spatial discretization methods.
	A relevant study concerning spherical Landau models is the calculation of multi-component lipid vesicles on a
	spherical surface \cite{Luo2018SPH_mini}.
	This study discretizes the modified
	Landau-Ginzburg free energy using spherical harmonics and directly obtains stationary states using the Broyden-Fletcher-Goldfarb-Shanno algorithm.
	Furthermore, a similar idea has been shown to be more efficient than the
	second category approach for finding stationary states of single- and
	multi-component phase field models in Euclidean space \cite{JiangNBPG,2021JiangBPG}.
	Motivated by this, our work aims to efficiently compute stationary states of the spherical LB model using the third type of numerical approach.
	\par
	This work has two main contributions.
	The first contribution is to develop a set of optimization algorithms to
	directly compute stationary states of discretized spherical LB free energy based on the spherical harmonics discretization.
	The optimization approaches include the accelerated adaptive-Bregman
	proximal gradient (AA-BPG), Nesterov, adaptive Nesterov (ANesterov),
	adaptive gradient descent (AGD) and adaptive nonlinear conjugate
	(ACG) methods.
	Theoretically, we give the convergence properties of these algorithms
	for the spherical LB model.
	Besides the efficient algorithms, good initial values can greatly \TR{speed the convergence} to stationary structures.
	Inappropriate initial values could lead to slow
	convergence rates, disordered states, or divergence.
	Our second contribution is to propose a principal mode analysis (PMA) method for good initial estimations to obtain the desired stationary configurations.
	The PMA method utilizes several principal spherical harmonics to capture the primary characteristics of the equilibrium structures.
	Furthermore, this approach reveals the relationship between the optimal sphere radius and the principal mode.
	Numerical results demonstrate that the PMA
	method can effectively estimate good initial values to speed up the process of finding stationary
	states.
	\par
	The rest of this article is organized as follows. 
	In Sect. \ref{sec:Model}, we introduce spherical LB model and spherical
	harmonics and then \FR{establish} the discretization formulation.
	In Sect. \ref{sec:Numer}, we develop optimization approaches to compute stationary states of spherical LB model, 
	and also present the procedures of the PMA method for estimating good initial values.
	In Sect. \ref{sec:Analysis}, we take some numerical experiments to
	demonstrate the power of developed approaches.
	Finally, we present concluding remarks and further developments in Sect. \ref{sec:Conc}.

%% file: Part_2_Problem_formulation.tex
	\section{Problem formulation}
	\label{sec:Model}
	We introduce the notations used throughout the paper.
	Let $\bbS^{2}:=\left( R\cos\theta \cos\phi, R\cos\theta \sin\phi,R\sin\theta
	\right)$ be the 2-dimensional spherical surface of radius $R$, where $\theta \in [0,\pi]$, $\phi \in[0,2\pi]$ are latitude and longitude angles, respectively.
	\tR{Spherical harmonics are solutions of the Laplace’s Equation in spherical coordinates. 
		Conveniently,
		spherical harmonics are constructed using Associated Legendre Polynomials \cite{2012Spherical}.
		A spherical harmonic function of degree $\ell$ and order $m$ is written as $Y^{m}_{\ell}(\theta,\phi)$.}
	We denote the Hilbert space $L^{2}(\bbS^{2})$ by $L^{2}$ in short, which includes all integrable periodic functions defined on a spherical surface.
	The inner product in $L^{2}$ is denoted by $\langle \cdot, \cdot \rangle$. 
	Let $\|\cdot\|$ and $\|\cdot\|_{\infty}$ be the $L^{2}$- and $L^{\infty}$-norm, respectively.
	More details about the Sobolev space and spherical harmonics can refer to \cite{2012Spherical}.
	
	\subsection{Spherical LB model}
	The spherical LB energy functional has the form of
	\begin{equation}
		\begin{aligned}
			E[\varphi]=\frac{1}{|\bbS^{2}|}\int_{\bbS^{2}}
			\left\{
			\frac{\xi^{2}}{2}[(1+\Delta)\varphi]^{2}+\frac{\epsilon}{2}\varphi^{2}-\frac{\lambda}{3!}\varphi^{3}+\frac{1}{4!}\varphi^{4}
			\right\} \; d\sigma,
		\end{aligned}
		\label{eq:LB_model}
	\end{equation}
	where $\varphi(R\cos\theta \cos\phi, R\cos\theta \sin\phi,R\sin\theta)$ is the order parameter describing the order
	degree in the system, while $\bbS^{2}$ is the spherical surface \fR{and $d\sigma$ is the infinitesimal element of $\mathcal{S}^2$}.
	The physical parameter $\xi$ corresponds to the bare correlation length, $\epsilon$ is a
	temperature-like variable, and $\lambda$ is a phenomenological parameter.
	The differential term represents the interaction potential, while the polynomial term represents the internal energy.
	Since the order parameter is the deviation from the average
	density \cite{1975Brazovskii}, the mass conservation holds
	\begin{equation}
		\begin{aligned}
			\frac{1}{|\bbS^{2}|} \int_{\mathbb{S}^{2}} \varphi \; d\sigma=0.
		\end{aligned}
		\label{eq:constrained_condition}
	\end{equation}

	\subsection{Spherical harmonics discretization}
	\label{sec:SPH}
	The purpose of this section is to give the mathematical formulation for searching stationary states of \eqref{eq:LB_model} subjected to the mass conservation constraint \eqref{eq:constrained_condition}.
	Here we use the spherical harmonic pseudo-spectral method to discretize the free energy functional \eqref{eq:LB_model}.
	\par
	We rewrite the order parameter function  $\varphi(\theta,\phi):=\varphi(R\cos\theta \cos\phi, R\cos\theta \sin\phi,R\sin\theta)$.
	Assume $\varphi(\theta,\phi)\in
	L^{2}(\bbS^{2})$, then $\varphi$ can be expanded by
	\begin{equation}
		\begin{aligned}
			\varphi(\theta,\phi)
			=\sum\limits^{\infty}_{\ell=0}\sum\limits_{m=-\ell}^{\ell}
			\hat{\varphi}_{\ell,m}Y^{m}_{\ell}(\theta,\phi).
		\end{aligned}
		\label{eq:order_parameter1}
	\end{equation}
	The spherical coefficients $\hat{\varphi}_{\ell,m}$ are calculated by
	\begin{equation}
		\begin{aligned}
			\hat{\varphi}_{\ell,m}=\int^{2\pi}_{0} \int^{\pi}_{0}\varphi(\theta,\phi) (Y^{m}_{\ell}(\theta,\phi))^{*} \sin \theta \; d\theta d\phi,
		\end{aligned}
		\label{eq:order_parameter_SH1}
	\end{equation}
	where $(Y^{m}_{\ell})^{*}=(-1)^{m}Y^{-m}_{\ell}$ is the conjugate function of $Y^{m}_{\ell}$. 
	The inner product $\langle \cdot,\cdot \rangle $ on a unit sphere is defined by
	\begin{equation}
		\begin{aligned}
			\langle f,g\rangle	=\int_{\bbS^{2}} f \bar{g}\; d\sigma
			=\int^{\pi}_{0}\int_{0}^{2\pi} f(\theta,\phi)\;
			\bar{g}(\theta,\phi)\;  \sin\theta \;
			d\phi\; d\theta,
		\end{aligned}
		\label{eq:L2_inner_product}
	\end{equation}
	where $\bar{g}$ denotes the complex conjugate function of $g$.
	We define the norms by $\|f\|=\langle f,f \rangle ^{1/2}$ and $\| f\|_{\infty}=\langle \sup f(\theta,\phi),1\rangle=\fR{\mathop{\max}\limits_{\ell,m}|\hat{f}_{\ell,m}|}$, \fR{where $\hat{f}_{\ell,m}$ is the corresponding spherical coefficient of $f$.} 
	Spherical harmonics have three important properties:
	\begin{enumerate}[(1)]
		\item Orthogonality
		\begin{equation}
			\begin{aligned}
				\langle Y^{m}_{\ell}, Y^{m^{'}}_{\ell^{'}}\rangle=\delta_{\ell\ell^{'}}\delta_{mm^{'}}.
			\end{aligned}
			\label{eq:orthogonality}
		\end{equation}
		\item  Completeness \\
		Let \tR{$\tilde{\varphi}_{N}$} be
		\begin{equation}
			\begin{aligned}
				\tR{\tilde{\varphi}_{N}(\theta,\phi)}=\sum^{N}_{\ell=0}\sum^{\ell}_{m=-\ell}
				\hat{\varphi}_{\ell,m} Y^{m}_{\ell}(\theta,\phi),
				\label{eq:finite_SHT}
			\end{aligned}
		\end{equation}	
		then \tR{$\tilde{\varphi}_{N}$} converges to $\varphi$ in the sense of $L^{2}-$ norm when $N\rightarrow \infty$. 
		\item Intrinsic definition \\
		The normalized spherical harmonics are  eigenfunctions of the spherical Laplacian operator $\Delta$ on the sphere with radius $R$
		\begin{equation}
			\begin{aligned}
				\Delta Y^{m}_{\ell}(\theta,\phi)=-\frac{\ell(\ell+1)}{R^{2}} Y^{m}_{\ell}(\theta,\phi).
			\end{aligned}
		\end{equation}
		
	\end{enumerate}
	\TR{\begin{remark}
			For $N$ in Eq.\;\eqref{eq:finite_SHT} sufficiently large, the error estimation of spatial discretization between   $\tilde{\varphi}_{N}(\theta,\phi)$ and $\varphi(\theta,\phi)$ in the form of Eq.\;\eqref{eq:order_parameter1} holds
			\begin{equation}
				\|\tilde{\varphi}_{N}(\theta,\phi)-\varphi(\theta,\phi)\|_{H^{t}(\mathcal{S}^2)}\leq (N+\frac{3}{2})^{-(s-t)}\|\varphi\|_{H^{s}(\mathcal{S}^2)}, \nonumber
			\end{equation}
			where $0 \leq t \leq s$, $H^{t}$ is the Sobolev space on a sphere and $t$ relates to the smoothness of the function.
			The proof and more details about Sobolev space on the sphere please refer to \cite{2012Spherical}.
	\end{remark}}
	
	According to the definition of spherical harmonics,
	the multiplication of two spherical harmonics holds
	\begin{equation}
		\begin{aligned}
			Y^{m_{1}}_{\ell_{1}} Y^{m_{2}}_{\ell_{2}} 
			=\sum\limits_{\ell \geq  0}\sum\limits_{\vert m \vert
				\leq \ell}\sqrt{\frac{(2\ell_{1}+1)(2\ell_{2}+1)(2\ell+1)}{4\pi}}
			\begin{pmatrix} \ell_{1} & \ell_{2} & \ell \\ 0 & 0 & 0 \end{pmatrix}
			\begin{pmatrix} \ell_{1} & \ell_{2} & \ell \\ m_{1} & m_{2} & m \end{pmatrix}
			Y^{m}_{\ell},
			\nonumber
		\end{aligned}
		\label{eq:SPH_2inner}
	\end{equation}
	where $\begin{pmatrix} \ell_{1} & \ell_{2} & \ell \\ m_{1} & m_{2} & m \end{pmatrix}$ is a Winger 3-j matrix \cite{Schulten1975wigner}.
	For convenience, let
	\begin{equation}
		\begin{aligned}
			C_{m_{1}m_{2}m}^{\ell_{1}\ell_{2}\ell}=\sqrt{\frac{(2\ell_{1}+1)(2\ell_{2}+1)(2\ell+1)}{4\pi}}
			\begin{pmatrix} \ell_{1} & \ell_{2} & \ell \\ 0 & 0 & 0 \end{pmatrix}
			\begin{pmatrix} \ell_{1} & \ell_{2} & \ell \\ m_{1} & m_{2} & m \end{pmatrix}.
		\end{aligned}
		\nonumber
	\end{equation}
	From the orthogonality \eqref{eq:orthogonality}, the triple and quadratic integrals hold
	\begin{equation}
		\begin{aligned}
			\langle Y^{m_{1}}_{\ell_{1}} Y^{m_{2}}_{\ell_{2}},
			Y^{m_{3}}_{\ell_{3}} \rangle
			&=\sqrt{\frac{(2\ell_{1}+1)(2\ell_{2}+1)(2\ell_{3}+1)}{4\pi}}
			\begin{pmatrix} \ell_{1} & \ell_{2} & \ell_{3} \\ 0 & 0 & 0 \end{pmatrix}
			\begin{pmatrix} \ell_{1} & \ell_{2} & \ell_{3} \\ m_{1} & m_{2} & m_{3} \end{pmatrix}
			\nonumber
			=  C^{\ell_{1}\ell_{2}\ell_{3}}_{m_{1}m_{2}m_{3}}, \\
			\langle Y^{m_{1}}_{\ell_{1}} Y^{m_{2}}_{\ell_{2}},
			Y^{m_{3}}_{\ell_{3}}  Y^{m_{4}}_{\ell_{4}} \rangle
			&=\sum\limits_{\ell,m}\sqrt{\frac{(2\ell_{1}+1)(2\ell_{2}+1)(2\ell_{3}+1)(2\ell_{4}+1)(2\ell+1)^{2}}{(4\pi)^{2}}}
			\nonumber
			\\
			& \times \begin{pmatrix} \ell_{1} & \ell_{2} & \ell \\ 0 & 0 & 0 \end{pmatrix}
			\begin{pmatrix} \ell_{1} & \ell_{2} & \ell \\ m_{1} & m_{2} & m \end{pmatrix}
			\begin{pmatrix} \ell_{3} & \ell_{4} & \ell \\ 0 & 0 & 0 \end{pmatrix}
			\begin{pmatrix} \ell_{3} & \ell_{4} & \ell \\ m_{3} & m_{4} & -m \end{pmatrix}
			\\
			&=\sum_{\ell,m}C_{m_{1}m_{1}m}^{\ell_{1}\ell_{2}\ell} C_{m_{3}m_{4}-m}^{\ell_{3}\ell_{4}\ell},
		\end{aligned}
		\label{eq:SPH_4inner}
	\end{equation}
	To make non-linear terms non-zero, from the properties of the Wigner 3-j matrix \cite{Dahlen_wigner_1998}, the indices of $C_{m_{1}m_{2}m_{3}}^{\ell_{1}\ell_{2}\ell_{3}}$ satisfy the following selection rules
	\begin{equation}
		\begin{aligned}
			& m_{1}+m_{2}+m_{3}=0,\; \ell_{1}+\ell_{2}+\ell_{3}=even,\; \vert \ell_{2}-\ell_{1} \vert \leq \ell_{3} \leq
			\ell_{1}+\ell_{2}.
		\end{aligned}
		\label{eq:condition_degree}
	\end{equation}
	$C_{m_{1}m_{2}m}^{\ell_{1}\ell_{2}\ell}$ and $C_{m_{3}m_{4}-m}^{\ell_{3}\ell_{4}\ell}$ obey the same rules. 	
	\par 
	From the above computation, the constrained spherical LB energy functional is discretized to the constrained finite dimensional optimization problem
	\begin{equation}
		\begin{aligned}
			\min\limits_{\hat{\varphi}_{\ell,m}\in X_{N} }\; E_{h}(\{\hat{\varphi}_{\ell,m}\}):=G_{h}(\{\hat{\varphi}_{\ell,m}\})+F_{h}(\{\hat{\varphi}_{\ell,m}\}), 
			\quad \mbox{s.t. $\hat{\varphi}_{0,0}=0$.}
		\end{aligned}
		\label{eq:D_LB_model}
	\end{equation}
	The feasible space $X_{N}$ satisfies
	\begin{equation}
		\begin{aligned}
			X_{N}:=\Big\{\{\hat{\varphi}_{\ell,m}\}_{0 \leq \ell\leq N, |m|\leq \ell}:\; \|\varphi\|^{2}=\sum\limits_{\ell,m} 
			|\hat{\varphi}^{m}_{\ell}|^{2} <\infty \Big\}.
		\end{aligned}
	\end{equation}
	The constraint condition in \eqref{eq:D_LB_model} is obtained by substituting Eq.\;\eqref{eq:finite_SHT} into Eq.\;\eqref{eq:constrained_condition}.	$G(\{\hat{\varphi}_{\ell,m}\})$ and $F(\{\hat{\varphi}_{\ell,m}\})$ have the form of
	\begin{equation}
		\begin{aligned}
			G_{h}(\{\hat{\varphi}_{\ell,m}\}) &=\frac{\xi^2}{2}\sum\limits_{\ell,m}
			(1-\frac{\ell(\ell+1)}{R^2})^2\hat{\varphi}_{\ell,m}\hat{\varphi}_{\ell,-m},
			\\
			F_{h}(\{\hat{\varphi}_{\ell,m}\}) &=\frac{\epsilon}{2}
			\sum\limits_{\ell,m}\hat{\varphi}_{\ell,m}\hat{\varphi}_{\ell,-m}
			-\frac{\lambda}{3!}
			\sum\limits_{\{\ell_{i},m_{i}\}_{i=1}^{3}}
			\TR{C^{\ell_{1}\ell_{2}\ell_{3}}_{m_{1}m_{2}m_{3}}}
			\hat{\varphi}_{\ell_{1},m_{1}}\hat{\varphi}_{\ell_{2},m_{2}}\hat{\varphi}_{\ell_{3},m_{3}}
			\\
			&+\frac{1}{4!}
			\sum\limits_{\{\ell_{i},m_{i}\}^{4}_{i=1},\ell,m} \TR{C^{\ell_{1}\ell_{2}\ell}_{m_{1}m_{2}m} C^{\ell_{3}\ell_{4}\ell}_{m_{3}m_{4}-m}}
			\hat{\varphi}_{\ell_{1},m_{1}}\hat{\varphi}_{\ell_{2},m_{2}}\hat{\varphi}_{\ell_{3},m_{3}}\hat{\varphi}_{\ell_{4},m_{4}} .
		\end{aligned}
		\label{eq:D_LB_model2}
	\end{equation}
	It is expensive to directly calculate the nonlinear term $F_{h}$ due to the
	convolutions of spherical harmonic coefficients.
	However, these convolutions are dot multiplications in physical space.
	Therefore, we can effectively calculate these convolutions via
	the discrete spherical harmonic transformation, implemented by the SHTns package \cite{shtns_Ishioka_2018}.

%% file: Part_3_Numerical_methods.tex
	\section{Numerical methods}
	\label{sec:Numer}
	
	In this section, we develop a series of efficient optimization methods
	to quickly find stationary states of \eqref{eq:D_LB_model}, and the PMA method to estimate good initial values.

	\subsection{Optimization methods}
	\label{subsec:Numer_LB}
	
	The spherical LB energy function
	$E_{h}(\{\hat{\varphi}_{\ell,m}\})=G_{h}(\{\hat{\varphi}_{\ell,m}\})+F_{h}(\{\hat{\varphi}_{\ell,m}\})$
	in \eqref{eq:D_LB_model} is continuously differentiable, where
	$\{\hat{\varphi}_{\ell,m}\}=\left(
	\hat{\varphi}_{0,0},\hat{\varphi}_{1,-1},\hat{\varphi}_{1,0},\hat{\varphi}_{1,1},\cdots,\hat{\varphi}_{N,N}\right)^{\mathrm{T}}$
	is the column vector of spherical coefficients.
	Let $\mathcal{S}=\{
	\{\hat{\varphi}_{\ell,m}\}:e_{0,0}^{\mathrm{T}}\{\hat{\varphi}_{\ell,m}\}=0\}$,
	where $e_{0,0}=(1,0,\cdots,0)^{\mathrm{T}}$,
	and $\delta_{\mathcal{S}}({\{\hat{\varphi}}_{\ell,m}\})=0$ if $\{\hat{\varphi}_{\ell,m}\}\in \mathcal{S}$, 
	otherwise $\delta_{\mathcal{S}}({\{\hat{\varphi}}_{\ell,m}\})=+\infty$.
	By choosing $\{\hat{\varphi}_{\ell,m}\} \in X_{N}\cap \mathcal{S}$,
	then the constrained problem \eqref{eq:D_LB_model} is reduced to a classical
	unconstrained nonconvex composite minimization problem, where $G_{h}$ is
	convex and $F_{h}$ is nonconvex.
	\par
	To solve such an optimization problem, numerical optimization methods can be employed.
	In present work, we develop the AGD, ACG, AA-BPG, Nesterov and ANesterov algorithms to the discretized spherical LB model $E_{h}(\{\hat{\varphi}_{\ell,m}\})=G_{h}(\{\hat{\varphi}_{\ell,m}\})+F_{h}(\{\hat{\varphi}_{\ell,m}\})$, as discussed in Sects. \ref{subsubsec:AGD_ACG} and \ref{subsubsec:Nesterov_BPG}.
	In the iterative process of these algorithms, appropriate step sizes need to be chosen to update the order parameter efficiently.
	The Nesterov method uses a fixed step size, while others employ a line search strategy to adaptively update the step sizes.
	
	\subsubsection{Line search}
	\label{subsec:line_search}
	
	In each step, the line search is initialized by a BB step \cite{BB_Barzilai_1988}
	\begin{equation}
		\begin{aligned}
			\alpha_{n}=\frac{\langle d_{n},\; d_{n} \rangle }{\langle d_{n},\; e_{n} \rangle},   \quad \mbox{or} \quad
			\alpha_{n}=\frac{\langle e_{n},\; d_{n} \rangle }{\langle e_{n},\; e_{n} \rangle},
			\label{eq:BB_step}
		\end{aligned}
	\end{equation}
	where \tR{$\langle,\rangle$ is the inner product defined by Eq.\;\eqref{eq:L2_inner_product},}
	\tR{ $d_{n}=\{\hat{\varphi}^{n}_{\ell,m}\}-\{\hat{\varphi}^{n-1}_{\ell,m}\}$}, 
	$e_{n}=\nabla E_{h}(\{\hat{\varphi}_{\ell,m}^{n}\})- \nabla E_{h}(\{\hat{\varphi}_{\ell,m}^{n-1}\})$. 
	$\nabla E_{h}(\{\hat{\varphi}_{\ell,m}^{n}\})$ denotes as the first derivative of $E_{h}$, i.e., $\nabla E_{h}=\nabla G_{h}+\nabla F_{h}$, where
	\begin{equation}
		\begin{aligned}
			\nabla G_{h}(\{\hat{\varphi}_{\ell,m}^{n}\})_{\ell,m}&=\xi^{2}\left(1-\frac{\ell(\ell+1)}{R^{2}}\right)\hat{\varphi}_{\ell,m}^{n}, \\
			\nabla F_{h}(\{\hat{\varphi}_{\ell,m}^{n}\})_{\ell,m}&=\epsilon \hat{\varphi}^{n}_{\ell,m}-\frac{\lambda}{2}\widehat{(\varphi^{n})^{2}}_{\ell,m}+\frac{1}{6}\widehat{(\varphi^{n})^{3}}_{\ell,m}.
		\end{aligned}
		\label{eq:D_LB_model_var}
	\end{equation}
	The non-convexity of $E_{h}$ could result in a negative step size.
	Thus we set $0<\alpha_{\min}\leq \alpha_{n}$ to avoid this.
	Algorithm \ref{alg:line-search} gives the procedures of the line search for
	finding the step sizes.
	\begin{algorithm}
		\caption{Line search for step size $\alpha_{n}$}
		\label{alg:line-search}
		\begin{algorithmic}[1]
			\REQUIRE Energy function $E_{h}$, \tR{iteration direction $p_{n}$}, the n-th $\{\hat{\varphi}_{\ell,m}^{n}\}$, $\rho\in(0,1)$ and $\alpha_{0}, \alpha_{\min}, \alpha_{\max}>0$ 
			\IF {$n=0$} 
			\STATE $\alpha_{n}=\alpha_{0}$
			\ELSE		
			\STATE Initialize $\alpha_{n}$ by Eq.\;\eqref{eq:BB_step} 
			\tR{\STATE  $\{\hat{\varphi}^{n+1}_{\ell,m}\}=\{\hat{\varphi}^{n}_{\ell,m}\}+\alpha_{n}p_{n}$}
			\fR{\STATE Compute $p_{n+1}$}
			\WHILE {$E_{h}(\{\hat{\varphi}_{\ell,m}^{n+1}\})>E_{h}(\{\hat{\varphi}_{\ell,m}^{n}\})$ \fR{or $p_{n+1}^{\mathrm{T}}\nabla E_{h}(\{\hat{\varphi}^{n+1}_{\ell,m}\})>0$} }
			\STATE  $\alpha_{n}=\rho \alpha_{n}$
			\STATE \tR{ $\{\hat{\varphi}^{n+1}_{\ell,m}\}=\{\hat{\varphi}^{n}_{\ell,m}\}+\alpha_{n}p_{n}$}
			\fR{\STATE Compute $E_{h}(\{\hat{\varphi}_{\ell,m}^{n+1}\})$, $\nabla E_{h}(\{\hat{\varphi}_{\ell,m}^{n+1}\})$, $p_{n+1}$}
			\fR{\IF {$\alpha_{n}\leq \alpha_{min}$}
				\STATE Break
				\ENDIF
			}
			\ENDWHILE
			\ENDIF
			\STATE \fR{Output} $\alpha_{n}=\max(\min(\alpha_{n},\alpha_{\max}),\alpha_{\min})$
		\end{algorithmic}  
	\end{algorithm}

	\subsubsection{AGD and ACG algorithms}
	\label{subsubsec:AGD_ACG}
	In this subsection, we introduce the AGD and ACG algorithms for searching stationary states of \eqref{eq:D_LB_model}, as outlined in Algorithm \ref{alg:optimization-SLB}.
	\fR{It is noted that we introduce a large enough number $P$ to bound the gradient error in Algorithm \ref{alg:optimization-SLB} in order to guarantee the convergence of the ACG method, such as $P=100$.}
	Furthermore, we prove the convergence properties of these algorithms, as demonstrated in Theorem \ref{theorem:AGD} and Theorem \ref{theorem:ACG}.
	
	\begin{algorithm} 
		\caption{ AGD and ACG methods}
		\label{alg:optimization-SLB}
		\begin{algorithmic}[1]
			\REQUIRE Initial values: $\{\hat{\varphi}_{\ell,m}^{0}\}$, $n=0$, $\alpha_{0}>0$, $\alpha>0$, $P>0$, $p_{0}=-\nabla E_{h}(\{\hat{\varphi}^{0}\})_{\ell,m}$, and $\tau>0$,  $n_{tol}\in \bbZ^{+}$
			\WHILE{$\|\nabla E_{h}(\{\hat{\varphi}_{\ell,m}^{n}\})\| \geq \tau$ }
			\IF{AGD method}
			\STATE Estimate $\alpha_{n}$,\;
			$\hat{\varphi}_{\ell,m}^{n+1}=\hat{\varphi}_{\ell,m}^{n}-\alpha_{n}\nabla E_{h}(\{\hat{\varphi}_{\ell,m}^{n}\})_{\ell,m}$ by Algorithm \ref{alg:line-search}
			\ENDIF
			\IF{ACG method}
			\STATE Estimate $\alpha_{n}$,\; $\hat{\varphi}_{\ell,m}^{n+1}=\hat{\varphi}_{\ell,m}^{n}+\alpha_{n}p_{n}$ by Algorithm \ref{alg:line-search} 
			\STATE
			$\beta_{n}=\|\left( \nabla E_{h}(\{\hat{\varphi}^{n+1}_{\ell,m}\})-\nabla E_{h}(\{\hat{\varphi}^{n}_{\ell,m}\})
			\right)^{\mathrm{T}}\nabla E_{h}(\{\hat{\varphi}^{n+1}_{\ell,m}\})  \|/\|\nabla E_{h}(\{\hat{\varphi}_{\ell,m}^{n}\})\|^{2}$ 
			\STATE
			$p_{n+1}=-\nabla E_{h}(\{\hat{\varphi}^{n+1}_{\ell,m}\})+\beta_{n}p_{n}$
			\STATE 
			if $\|p_{n+1}\|>P$ \fR{or $p_{n+1}^{\mathrm{T}}\nabla E_{h}(\{\hat{\varphi}^{n+1}\})>0$,} reset $p_{n+1}=-\nabla E_{h}(\{\hat{\varphi}^{n+1}_{\ell,m}\})$
			\ENDIF
			\STATE $n=n+1$
			\IF{$n>n_{tol}$}
			\STATE Break
			\ENDIF
			\ENDWHILE
		\end{algorithmic}  
	\end{algorithm}
	\TR{
		\begin{remark}
			We should point out that in the practical simulation, we set $\hat{\varphi}^{n+1}_{0,0}=0$ for all $n$ after updating $\hat{\varphi}^{n+1}_{\ell,m}$ via step 3 (or step 6) in Algorithm\; 2 to guarantee the mass conservation in AGD and ACG methods.
		\end{remark}
	}	
	\begin{theorem}
		\label{theorem:AGD}
		Consider the AGD scheme $\hat{\varphi}^{n+1}=\hat{\varphi}^{n}+\alpha_{n} p_{n}$ with $p_{n}:=-\nabla E_{h}(\hat{\varphi}^{n})$.
		The following assumptions hold
		\begin{enumerate}[(i)]
			\item $E_{h}$ is bounded below in $X_{N}$ and \TR{$\mathcal{L}:=\{\hat{\varphi}: E_{h}(\hat{\varphi})\leq E_{h}(\hat{\varphi}^{0})\}$}, where $\hat{\varphi}^{0}$ is the initial value of the iteration. \\
			\item $\nabla E_{h}$ is Lipschitz continuous in any open subset $\mathcal{N}$ of $\mathcal{L}$, i.e., there exists a positive constant $L$ such that
			\begin{equation}
				\begin{aligned}
					\|\nabla E_{h}(\psi)-\nabla E_{h}(\tilde{\psi})\|\leq L\|\psi-\tilde{\psi}\|, \quad\quad \mbox{$\psi, \tilde{\psi}\in \mathcal{N}$.}
				\end{aligned}
			\end{equation}			
			\item The step length $\alpha_{n}$ satisfies the Wolfe condition 
				\begin{equation}   
					\begin{aligned}
						&E_{h}(\hat{\varphi}^{n}+\alpha_{n}p_{n})\leq E_{h}(\hat{\varphi}^{n})+ c_{1}\alpha_{n}\nabla p_{n}^{\mathrm{T}}  E_{h}(\hat{\varphi}^{n}) , \\
						& p_{n}^{\mathrm{T}}\nabla E_{h}(\hat{\varphi}^{n}+\alpha_{n}\tR{p_{n}})\geq c_{2}p_{n}^{\mathrm{T}}
						\tR{\nabla} E_{h}(\hat{\varphi}^{n}),
					\end{aligned}
					\label{eq:Wolfe_condition}
				\end{equation}
		\end{enumerate} 
		where $0<c_{1}<c_{2}<1$. Then we have   
		\begin{equation}
			\begin{aligned}
				\displaystyle\lim_{n\rightarrow \infty}\inf \|\nabla  E_{h}(\hat{\varphi}^{n})\| =0.
			\end{aligned}
			\label{eq:AGD_result}
		\end{equation}
		\begin{proof}
			According to \eqref{eq:Wolfe_condition}, we have
			\begin{equation}
				\begin{aligned}
					p_{n}^{\mathrm{T}}\left(\nabla E_{h}(\hat{\varphi}^{n}+\alpha_{n}p_{n})-\nabla E_{h}(\hat{\varphi}^{n})\right)  \geq (c_{2}-1)p_{n}^\mathrm{T}\nabla E_{h}(\hat{\varphi}^{n}) .
				\end{aligned}
				\label{eq:proof_1}
			\end{equation}
			Due to the Lipschitz continuity of $\nabla E_{h}$ and $c_{2}<1$, it becomes
			\begin{equation}
				\begin{aligned}
					(c_{2}-1)p_{n}^{\mathrm{T}} \nabla E_{h}(\hat{\varphi}^{n})  \leq \alpha_{n}\cdot L \|p_{n}\|^{2}.
				\end{aligned}
			\end{equation}
			Thus the step length at each step fulfils
			\begin{equation}
				\begin{aligned}
					\alpha_{n}\geq \frac{c_{2}-1}{L}\cdot \frac{p_{n}^{\mathrm{T}}\nabla E_{h}(\hat{\varphi}^{n})} {\|p_{n}\|^{2}} >0.
				\end{aligned}
				\label{eq:AGD_alpha}
			\end{equation}
			By substituting \eqref{eq:AGD_alpha} into \eqref{eq:Wolfe_condition}, we obtain
			\begin{equation}
				\begin{aligned}
					E_{h}(\hat{\varphi}^{n+1})\leq E_{h}(\hat{\varphi}^{n})-\frac{c_{1}(1-c_{2})}{L}\frac{\|p_{n}^{\mathrm{T}}\nabla E_{h}(\hat{\varphi}^{n}) \|^{2}}{\|p_{n}\|^{2}}.
				\end{aligned}
				\label{eq:AGD_ineq}
			\end{equation}
			With $p_{n}=-\nabla E_{h}(\hat{\varphi}^{n})$, summing the expression over all $n\geq 0$, Eq.\;\eqref{eq:AGD_ineq} leads to
			\begin{equation}
				\begin{aligned}
					E_{h}(\hat{\varphi}^{n+1})\leq E_{h}(\hat{\varphi}^{0})-\sum^{\infty}_{n=0}\frac{c_{1}(1-c_{2})}{L} \|\nabla E_{h}(\hat{\varphi}^{n})\|^{2}.
				\end{aligned}
			\end{equation}
			Since $E_{h}$ is bounded, we have
			\begin{equation}
				\begin{aligned}
					\sum^{\infty}_{n=0}\|\nabla E_{h}(\hat{\varphi}\tR{^{n}})\|^{2} < \infty.
				\end{aligned}
			\end{equation}
			This completes the proof of the convergence of the AGD algorithm for the
			spherical LB model.
		\end{proof}
		
	\end{theorem}
	
	\begin{theorem}
		\label{theorem:ACG}
		Consider the ACG scheme
		\begin{equation}
			\begin{aligned}
				&\hat{\varphi}^{n+1}=\hat{\varphi}^{n}+\alpha_{n} p_{n},\\
				&p_{n+1}=-\nabla E_{h}(\hat{\varphi}^{n+1})+\beta_{n}p_{n}, \\
				&\beta_{n}= \min\left( \left|
				\left(\nabla E_{h}(\hat{\varphi}^{n+1})-\nabla
				E_{h}(\hat{\varphi}^{n}) \right)^{\mathrm{T}} \nabla E_{h}(\hat{\varphi}^{n+1}) / \|\nabla
				E_{h}(\hat{\varphi}^{n})\|^{2} \right|
				, \|\nabla E_{h}(\hat{\varphi}^{n+1})\|^{2} / \|\nabla E_{h}(\hat{\varphi}^{n})\|^{2} \right).
			\end{aligned}
			\label{eq:ACG_iter}
		\end{equation}
		The following assumptions hold
		\begin{enumerate}[(i)]
			\item $E_{h}$ is bounded below in $X_{N}$ and
			\TR{ $\mathcal{L}:=\{\hat{\varphi}: E_{h}(\hat{\varphi})\leq E_{h}(\hat{\varphi}^{0})\}$} is bounded, where $\hat{\varphi}^{0}$ is the initial value of the iteration. \\
			\item $\nabla E_{h}$ is Lipschitz continuous
			differentiable in some open neighbourhood $\mathcal{N}$ of $\mathcal{L}$. \\
			\item The step length $\alpha_{n}$ satisfies the strong Wolfe condition 
				\begin{equation}    
					\begin{aligned}
						& E_{h}(\hat{\varphi}^{n}+\alpha_{n}p_{n})\leq E_{h}(\hat{\varphi}^{n})+ c_{1}\alpha_{n}p_{n}^{\mathrm{T}}\nabla E_{h}(\hat{\varphi}^{n}) , \\
						&|\nabla p_{n}^{\mathrm{T}} E_{h}(\hat{\varphi}^{n}+\alpha_{n}\tR{p_{n}})| \leq -c_{2}p_{n}^{\mathrm{T}}\tR{\nabla} E_{h}(\hat{\varphi}^{n}) ,
					\end{aligned}
					\label{eq:Wolfe_condition2}
				\end{equation}
		\end{enumerate} 			
		where $0<c_{1}<c_{2}< 1/2$. Then we have 
		\begin{equation}
			\begin{aligned}
				\displaystyle\lim_{n\rightarrow \infty}\inf \|\nabla E_{h}(\hat{\varphi}^{n})\| =0.
			\end{aligned}
			\label{eq:ACG_result}
		\end{equation}
		\begin{proof}
			\par
			We first prove that $p_{n}$ is a descent direction \fR{and there exists  $\alpha_{n}\in(0,\alpha_{n}^{\max})$ satisfying the strong Wolfe condition, i.e.}
			\begin{subequations}
				\begin{align}
						-\frac{1}{1-c_{2}} \leq 
						\frac{ p_{n}^{\mathrm{T}}\nabla E_{h}(\hat{\varphi}^{n})}{\|\nabla E_{h}(\hat{\varphi}^{n})\|^{2}} \leq \frac{2c_{2}-1}{1-c_{2}}, \quad \mbox {for $n=0,1,\cdots$},
						\label{eq:ACG_proof1_1}
					\end{align}
					\begin{align}
						\fR{\alpha_{n}^{\max}=\min\left\{ \frac{2(1-c_{1})}{L_{n}(1-c_{2})}, \frac{2(1-c_{1})\|\nabla E_{h}(\hat{\varphi}^{n})\|^{2}}{L_{n}(1-c_{2})P^{2}} ,1\right\}.}
					\label{eq:ACG_proof1_2}
				\end{align}
			\end{subequations}
			This expression can be proved by induction.  
			The procedure can be split into three steps.
			\par
			Step 1: We show that $p_{0}$ satisfies \eqref{eq:ACG_proof1_1}.
			For $n=0$, the middle term is $-1$ \fR{since $p_{0}=\nabla E(\hat{\varphi}^{0})$}. 
			In fact, according to $0<c_{2}<1$, we have
			\begin{equation}
				\begin{aligned}
					-1<\frac{2c_{2}-1}{1-c_{2}}<0.
				\end{aligned}
			\end{equation}
			Thus Eq.\;\eqref{eq:ACG_proof1_1} is valid for $n=0$.
			\par
			\fR{Step 2: We show there exists $\alpha_{0}$ in the form of Eq.\;\eqref{eq:ACG_proof1_2} satisfying the strong Wolfe condition Eq.\,\eqref{eq:Wolfe_condition2}.
				Generally, we denote $j=0$.
				From the Taylor expansion, we obtain
				\begin{equation}
					\begin{aligned}
						E_{h}(\hat{\varphi}^{j}+\alpha_{j}p_{j})
						&=E_{h}(\hat{\varphi}^{j})+\alpha_{j}\langle \nabla E_{h}(\hat{\varphi}^{j}),p_{j} \rangle+\frac{\alpha_{n}^{2}}{2}\langle \nabla^{2} E_{h}(\zeta^{j})p_{j},p_{j} \rangle,
						\\
						&=E_{h}(\hat{\varphi}^{j})+c_{1}\alpha_{j}\langle \nabla E_{h}(\hat{\varphi}^{j}),p_{j} \rangle+
						(1-c_{1})\alpha_{j}\langle \nabla E_{h}(\hat{\varphi}^{j}),p_{j} \rangle
						+\frac{\alpha_{j}^{2}}{2}\langle \nabla^{2} E_{h}(\zeta^{j})p_{j},p_{j} \rangle,
					\end{aligned}
					\label{eq:Wolfe_condition_proof1}
				\end{equation}
				where \fR{$\zeta^{j}\in V_{j}:=\{ \hat{\varphi}^{j}+\alpha p_{j}, \alpha \in (0,1)\}$.}
				With $\langle \nabla E_{h}(\hat{\varphi}^{j}),p_{j}\rangle <0$ and $0<c_{1}<1$, 
				$\alpha_{j}$ satisfies
				\begin{equation}
					\begin{aligned}
						\alpha_{j}\leq \frac{-2(1-c_{1})\langle \nabla E_{h}(\hat{\varphi}^{j}),p_{j} \rangle} {L_{j}\|p_{j}\|^{2}},
					\end{aligned}
					\label{eq:ACG_bound}
				\end{equation}
				where $L_{j}=\max\{\|\nabla^{2}E_{h}(x)\|:x\in V_{j}\}$.
				Note that for the ACG scheme in Algorithm \ref{alg:optimization-SLB}, we restrict $p_{j}=-\nabla E_{h}(\hat{\varphi}^{j})$ when $\|p_{j}\|>P$, where $P$ is a finite positive constant.
				Thus we have 
				\begin{equation}
					\begin{aligned}
						\alpha_{j}^{\max}=\min\left\{ \frac{-2(1-c_{1})\langle \nabla E_{h}(\hat{\varphi}^{j}),p_{j} \rangle} {L_{j}\|\nabla E_{h}(\hat{\varphi}^{j})\|^{2}},
						\frac{ -2(1-c_{1}) \langle \nabla E_{h}(\hat{\varphi}^{j}),p_{j} \rangle} {L_{j}P^{2}},1
						\right\}.
					\end{aligned}
					\label{eq:ACG_bound2}
				\end{equation}
				Thus Eq.\;\eqref{eq:ACG_proof1_2} is valid for $n=0$.
				This also implies
				\begin{equation*}
					(1-c_{1})\alpha_{j}\langle \nabla E_{h}(\hat{\varphi}^{j}),p_{j} \rangle
					+\frac{\alpha_{j}^{2}}{2}\langle \nabla^{2} E_{h}(\zeta^{j})p_{j},p_{j} \rangle \leq 0.
				\end{equation*}
				Therefore, $\alpha_{0}$ satisfies the first inequation in Wolfe condition \eqref{eq:Wolfe_condition2}.
				\par
				Next, we prove that the $\alpha_{0}$ also satisfies the second inequation in \eqref{eq:Wolfe_condition2}. 
				\fR{Since $E_{h}(\hat{\varphi}^{j}+\alpha p_{j})$ is bounded below and $E_{h}(\hat{\varphi}^{j})+\alpha c_{1}p_{j}^{\mathrm{T}}\nabla E_{h}(\hat{\varphi}^{j}) $ is unbounded below,
					they have at least one intersection point.
					Let $\alpha_{j}$ be the smallest intersection point, that is
					\begin{equation}
						\begin{aligned}					E_{h}(\hat{\varphi}^{j}+\alpha_{j}p_{j})			=E_{h}(\hat{\varphi}^{j})+\alpha_{j} c_{1}p_{j}^{\mathrm{T}}\nabla E_{h}(\hat{\varphi}^{j}) .
						\end{aligned}
						\label{eq:middle_equation2}
					\end{equation}
					According to the Taylor expansion,
					there exists \fR{$\zeta_{1}^{j}= \hat{\varphi}^{j}+\bar{\alpha} p_{j}$, $\bar{\alpha} \in (0,\alpha_{j})$ such that}
					\begin{equation}
						\begin{aligned}
							E_{h}(\hat{\varphi}^{j}+\alpha_{j}p_{j})
							=E_{h}(\hat{\varphi}^{j})+\alpha_{j} p_{j}^{\mathrm{T}}\nabla E_{h}(\zeta_{1}^{j}) .
						\end{aligned}
						\label{eq:middle_equation1}
				\end{equation}}
				By combining \eqref{eq:middle_equation1} and \eqref{eq:middle_equation2}, we have
				\begin{equation}
					\begin{aligned}
						p_{j}^{\mathrm{T}}\nabla E_{h}(\zeta_{1}^{j}) = c_{1}p_{j}^{\mathrm{T}}\nabla E_{h}(\hat{\varphi}^{j}) .\nonumber
					\end{aligned}
				\end{equation}
				From $\langle \nabla E_{h}(\hat{\varphi}^{j}),p_{j}\rangle <0$ and $0<c_{1}<c_{2}<1$, we have 
				\begin{equation}
					\begin{aligned}
						&p_{n}^{\mathrm{T}} \nabla E_{h}(\zeta_{1}^{j})  > c_{2}p_{j}^{\mathrm{T}}\nabla E_{h}(\hat{\varphi}^{j}) , \\
						& |p_{j}^{\mathrm{T}}\nabla E_{h}(\zeta_{1}^{j}) |<-c_{2}p_{j}^{\mathrm{T}}\nabla E_{h}(\hat{\varphi}^{j}).
					\end{aligned}
				\end{equation}		
				In summary, the step size $p_{0}$ satisfies the Wolfe condition or strong Wolfe condition.
				Since $p_{0}$ satisfies \eqref{eq:ACG_proof1_1}, we have
				\begin{equation*}
					-\frac{1}{1-c_{2}} \leq 
					\frac{p_{j}^{\mathrm{T}}\nabla E_{h}(\hat{\varphi}^{j}) }{\|\nabla E_{h}(\hat{\varphi}^{j})\|^{2}} \leq \frac{2c_{2}-1}{1-c_{2}}.
				\end{equation*}
				Thus \eqref{eq:ACG_bound2} becomes
				\begin{equation}
					\begin{aligned}
						\alpha_{j}^{\max}=\min\left\{ \frac{2(1-c_{1})}{L_{j}(1-c_{2})}, \frac{2(1-c_{1})\|\nabla E_{h}(\hat{\varphi}^{j})\|^{2}}{L_{j}(1-c_{2})P^{2}} ,1\right\}.
					\end{aligned}
					\label{eq:ACG_bound3}
				\end{equation}   
				To be end, we have shown that Eq.\;\eqref{eq:ACG_proof1_2} is valid for $n=0$.}
			\par
			\fR{
				Step 3: We show that Eqs.\;\eqref{eq:ACG_proof1_1}-\eqref{eq:ACG_proof1_2} also hold for $n>0$.}
			\fR{From \eqref{eq:ACG_iter}, let $n=0$, then the $n+1$ step holds}
			\begin{equation}
				\begin{aligned}
					\frac{ p_{n+1}^{\mathrm{T}} \nabla E_{h}(\hat{\varphi}^{n+1}) } {\| \nabla E_{h}(\hat{\varphi}^{n+1}) \|^{2}} =-1 +\beta_{n} \frac{ p_{n}^{\mathrm{T}} \nabla E_{h}(\hat{\varphi}^{n+1}) } {\| \nabla E_{h}(\hat{\varphi}^{n+1}) \|^{2}}.
				\end{aligned}
			\end{equation}
			\fR{Since $\alpha_{0}$ satisfies \eqref{eq:Wolfe_condition2}} and $|\beta_{n}|\leq \|\nabla E_{h}(\hat{\varphi}^{n+1})\|^{2} / \|\nabla E_{h}(\hat{\varphi}^{n})\|^{2}$, the above equation is reduced to 
			\begin{equation}
				-1+c_{2}\frac{ p_{n}^{\mathrm{T}} \nabla E_{h}(\hat{\varphi}^{n}) } {\| \nabla E_{h}(\hat{\varphi}^{n}) \|^{2}}
				\leq
				\frac{ p_{n+1}^{\mathrm{T}} \nabla E_{h}(\hat{\varphi}^{n+1}) } {\| \nabla E_{h}(\hat{\varphi}^{n+1}) \|^{2}}
				\leq 
				-1-c_{2}\frac{ p_{n}^{\mathrm{T}} \nabla E_{h}(\hat{\varphi}^{n})
				} {\| \nabla E_{h}(\hat{\varphi}^{n}) \|^{2}}.
			\end{equation}
			Since $p_{n}$ satisfies \eqref{eq:ACG_proof1_1}, then we have
			\begin{equation}
				\begin{aligned}
					-1-\frac{c_{2}}{1-c_{2}} \leq 
					\frac{p_{n+1}^{\mathrm{T}} \nabla E_{h}(\hat{\varphi}^{n+1}) }{\|\nabla E_{h}(\hat{\varphi}^{n+1})\|^{2}} \leq -1+\frac{c_{2}}{1-c_{2}}, \quad \mbox {for $n=0,1,\cdots$}.
				\end{aligned}
				\label{eq:ACG_proof2}
			\end{equation}
			This equation implies that \eqref{eq:ACG_proof1_1} holds for all $n=1$, \fR{i.e. $p_{1}$ is the descent direction.}
			\fR{Then repeat the Step 2, we can show that there exists $\alpha_{1}$ in the form of \eqref{eq:ACG_proof1_2} satisfying the strong Wolfe condition \eqref{eq:Wolfe_condition2}.
				By induction, we show that $p_{i}$ and $\alpha_{i}$ for $i\geq 2$ satisfy \eqref{eq:ACG_proof1_1}-\eqref{eq:ACG_proof1_2}.}		
			
			\par
			Finally, we prove the global convergence \eqref{eq:ACG_result} by contradiction. 
			Assume that there exists a positive constant $\kappa$ such that
			\begin{equation}
				\begin{aligned}
					\|\nabla E_{h}(\hat{\varphi}^{n})\| \geq \kappa,
				\end{aligned}
				\label{eq:ACG_hypothesis}
			\end{equation}
			\TR{for all $n$ sufficiently large.}
			\TR{
				From \eqref{eq:ACG_proof1_1} and $0<c_{2}\leq 1/2$, we can obtain 
				\begin{equation}
					\|p_{n}^{\mathrm{T}} \nabla E_{h}(\hat{\varphi}^{n})\|^{2}\geq \left(\frac{2c_{2}-1}{1-c_{2}}\right)^{2} \|\nabla E_{h}(\hat{\varphi}^{n})\|^{4}.
					\nonumber
				\end{equation}
				We recall Eq.\;\eqref{eq:AGD_ineq}, holding
				\begin{equation}
					\begin{aligned}
						E_{h}(\hat{\varphi}^{n+1})\leq E_{h}(\hat{\varphi}^{n})-\frac{c_{1}(1-c_{2})}{L}\frac{\|p_{n}^{\mathrm{T}}\nabla E_{h}(\hat{\varphi}^{n}) \|^{2}}{\|p_{n}\|^{2}}. \nonumber
					\end{aligned}
				\end{equation}
				Then we have
				\begin{equation*}
					E_{h}(\hat{\varphi}^{n+1})\leq E_{h}(\hat{\varphi}^{n})-\frac{c_{1}(1-c_{2})^2}{L(1-c_{2})}\frac{\| \nabla E_{h}(\hat{\varphi}^{n})\|^{4}} {\|p_{n}\|^{2}}.
				\end{equation*}
				Since $E_{h}$ is bounded below, summing the expression over all $n\geq 0$, }		
			we conclude
			\begin{equation}
				\begin{aligned}
					\sum^{\infty}_{n=0} \frac{\| \nabla E_{h}(\hat{\varphi}^{n})\|^{4}} {\|p_{n}\|^{2}} < \infty .
				\end{aligned}
				\label{eq:ACG_proof2}
			\end{equation}
			From Eqs.\;\eqref{eq:Wolfe_condition2} and \eqref{eq:ACG_proof1_1}, we obtain
			\begin{equation}
				\begin{aligned}
					|p_{n-1}^{\mathrm{T}}\nabla E_{h}(\hat{\varphi}^{n}) |\leq -c_{2} p_{n-1}^{\mathrm{T}} \nabla E_{h}(\hat{\varphi}^{n-1})  \leq \frac{c_{2}}{1-c_{2}}\|\nabla E_{h}(\hat{\varphi}^{n-1})\|^{2}.
				\end{aligned}
			\end{equation}
			Therefore, we obtain
			\begin{equation}
				\begin{aligned}
					\|p_{n}\|^{2}=\|-\nabla E_{h}(\hat{\varphi}^{n})+\beta_{n-1}p_{n-1} \|^{2}
					& \leq \|\nabla E_{h}(\hat{\varphi}^{n})\|^{2} +2|\beta_{n-1}||p_{n-1}^{\mathrm{T}}\nabla E_{h}(\hat{\varphi}^{n}) |+\beta_{n-1}^{2}\|p_{n-1}\|^{2} 
					\\
					&\leq \|\nabla E_{h}(\hat{\varphi}^{n})\|^{2} +\frac{2c_{2}}{1-c_{2}}|\beta_{n-1}|\|\nabla E_{h}(\hat{\varphi}^{n-1})\|^{2}+\beta_{n-1}^{2}\|p_{n-1}\|^{2} \\
					&\leq \frac{1+c_{2}}{1-c_{2}}\|\nabla E_{h}(\hat{\varphi}^{n})\|^{2}  +\frac{\|\nabla E_{h}(\hat{\varphi}^{n})\|^{2} }{\|\nabla E_{h}(\hat{\varphi}^{n-1})\|^{2}}\|p_{n-1}\|^{2},
				\end{aligned}
			\end{equation}
			where the third inequality is obtained from the definition $\beta_{n}$ in \eqref{eq:ACG_iter}.
			From $(1+c_{2})/(1-c_{2})\geq 1$ and applying this relation repeatedly to $p_{n-1},\cdots,p_{1}$, we have 
			\begin{equation}
				\begin{aligned}
					\|p_{n}\|^{2}\leq \frac{1+c_{2}}{1-c_{2}}   \|\nabla
					E_{h}(\hat{\varphi}^{n})\|^{4}\sum^{n}_{k=0}  \|\nabla
					E_{h}(\hat{\varphi}^{k})\|^{-2}.
				\end{aligned}
			\end{equation}
			The assumptions (i)-(ii) imply that there exists a  constant $\bar{\kappa}>0$ such that $\|\nabla E_{h}(\hat{\varphi}^{n})\|\leq \bar{\kappa}$ when $\hat{\varphi}^{n}\in \mathcal{L}$,
			then we have
			\begin{equation}
				\begin{aligned}
					\|p_{n}\|^{2}\leq \frac{1+c_{2}}{1-c_{2}}\frac{\bar{\kappa}^{4}}{\kappa^{2}}n.
				\end{aligned}
			\end{equation}
			Thus we obtain
			\begin{equation}
				\begin{aligned}				\sum^{\infty}_{n=1}\frac{1}{\|p_{n}\|^{2}} \geq \kappa_{1}\sum^{\infty}_{n=1}\frac{1}{n},
				\end{aligned}
			\end{equation}
			where $\kappa_{1}$ is a positive constant.
			Since the right side is $\infty$ when $n\rightarrow\infty$,
			then we have $\sum^{\infty}_{n=1}\frac{1}{\|p_{n}\|^{2}} > \infty$, which conflicts with \eqref{eq:ACG_proof2}.
			This implies that the hypothesis \eqref{eq:ACG_hypothesis} is not true.
			Thus the statement in Theorem \ref{theorem:ACG} is proved.
		\end{proof}   
	\end{theorem}
	\par 

	\fR{
		\begin{remark}
			For the AGD method, the iteration direction $p_{j}=-\nabla E_{h}(\hat{\varphi}^{j})$ is a descent direction.
			Thus from the Step 2 in Theorem \ref{theorem:ACG}, we can conclude that the maximum step size $\alpha_{j}^{\max}$ for AGD convergence theorem (Theorem \ref{theorem:AGD}) satisfies
			\begin{equation}
				\alpha_{j}^{\max}=\min\left\{ 2(1-c_{1})/L_{j},1\right\}, \quad \mbox{$j=0,1,2,\cdots$},
			\end{equation}
			where $L_{j}=\max\{\|\nabla^2 E_{h}(x)\|:x=\hat{\varphi}^{j}-\alpha \nabla E_{h}(\hat{\varphi}^{j}),\alpha\in(0,1)\}$.
		\end{remark}
	}	
	
	\subsubsection{AA-BPG and (A)Nesterov algorithms}
	\label{subsubsec:Nesterov_BPG}
	This subsection outlines the procedures and convergence properties of the AA-BPG, Nesterov and ANesterov algorithms for computing stationary states of \eqref{eq:D_LB_model}.
	The main \FR{challenges} in the Nesterov, ANesterov and AA-BPG algorithms \FR{lie} in efficiently solving the following equations 
	\begin{equation}
		\begin{aligned}
			\hat{\psi}_{\ell,m}^{n} &= \hat{\varphi}^{n}_{\ell,m}+w_{n}(\hat{\varphi}^{n}_{\ell,m}-\hat{\varphi}^{n-1}_{\ell,m}),\\
			\hat{\varphi}^{n+1}_{\ell,m} &= \arg \min_{\hat{\varphi}_{\ell,m}\in X_{N}}\left\{G_{h}(\{\hat{\varphi}_{\ell,m}\})+
			\langle \{\hat{\varphi}_{\ell,m}-\hat{\psi}^{n}_{\ell,m}\} , \nabla F_{h}(\{\hat{\psi}_{\ell,m}^{n}\})\rangle
			+\frac{1}{\alpha_{n}}D_{c}(\{\hat{\varphi}_{\ell,m}\},\{\hat{\psi}_{\ell,m}^{n}\}) \right\},\
			\\
			& \mbox{s.t. $\hat{\varphi}_{0,0}=0$}.
		\end{aligned}
		\label{eq:gen_psi}
	\end{equation}
	$D_{c}$ is a Bregman distance defined by a convex function $c(x)$
	\begin{equation}
		\begin{aligned}
			D_{c}(x,y)=c(x)-c(y)-\langle \nabla c(x),x-y\rangle,\quad 
			\mbox{$(x,y)\in\; dom\;c \times intdom\;c$,}
		\end{aligned}
		\label{eq:Bregman_distance}
	\end{equation}
	where $dom\;c=\{x: c(x)<\infty \}$ and $intdom\;c$ is the set consisting of all interior points of $dom\;c$.
	Two different convex functions $c(x)$ can be chosen as
	\begin{equation}
		\begin{aligned}
			&\mbox{M=2:}\quad\quad c(x)=\frac{1}{2}\|x\|^{2}, \\
			&\mbox{M=4:}\quad\quad c(x)=\frac{a}{4}\|x\|^{4}+\frac{b}{2}\|x\|^{2}+1, \quad \mbox{$a,b>0$,}
		\end{aligned}
		\label{eq:convex_function_c}
	\end{equation}
	where M denotes the highest order of the polynomial $c(x)$.
	According to the type of $c(x)$,
	we recall them as the AA-BPG-2 method and the AA-BPG-4 method, respectively. 
	For the AA-BPG-2 method, 
	\TR{the Euler-Lagrange equation corresponding to Eq.\;\eqref{eq:gen_psi} is formulated as}
	\begin{equation}
		\begin{aligned}
			&\alpha_{n} \nabla G_{h}(\{\hat{\varphi}^{n+1}_{\ell,m}\})+
			\alpha_{n}\nabla F_{h}(\{\hat{\psi}^{n}_{\ell,m}\})
			+(\{\hat{\varphi}^{n+1}_{\ell,m}\}-\{\hat{\psi}^{n}_{\ell,m}\})-\alpha_{n}\gamma_{n}e_{0,0}=0. \\
		\end{aligned}
	\end{equation}
	From the mass conservation $e_{0,0}^{\mathcal{T}}\{\hat{\varphi}_{\ell,m}\}=0$,
	we obtain the Lagrange multiplier $\gamma_{n}=\nabla F_{h}(\{\hat{\varphi}_{\ell,m}^{n}\})_{0,0}$.
	Due to $\nabla G_{h}(\{\hat{\varphi}_{\ell,m}^{n+1}\})_{\ell,m}=(1-  (\ell(\ell+1)/R^2) )\hat{\varphi}_{\ell,m}^{n+1}$,
	the above equation becomes
	\begin{equation}
		\begin{aligned}
			\{\hat{\varphi}^{n+1}_{\ell,m}\}=\left( \alpha_{n}D+\mathrm{I}\right)^{-1}
			\left( \{\hat{\psi}_{\ell,m}^{n}\}
			-\alpha_{n}\nabla F_{h}(\{\hat{\psi}_{\ell,m}^{n}\}) \right),
			\quad\quad \mbox{$\nabla F_{h}(\{\psi^{n}_{\ell,m}\})_{0,0}=0$,}
		\end{aligned}
		\label{eq:D_LB_model_G2_1}
	\end{equation}
	where $D$ is a diagonal matrix and the diagonal elements are $1-\ell(\ell+1)/R^{2}$.
	Eq.\;\eqref{eq:D_LB_model_G2_1} is the standard semi-implicit scheme (SIS). 
	Similarly, for the AA-BPG-4 method, we have the following \TR{Euler-Lagrange equation based on \eqref{eq:gen_psi}}
	\begin{equation}
		\begin{aligned}
			&\alpha_{n} \nabla G_{h}(\{\hat{\varphi}^{n+1}_{\ell,m}\})+
			\alpha_{n}\nabla F_{h}(\{\hat{\psi}^{n}_{\ell,m}\})
			+(a\|\{\hat{\varphi}_{\ell,m}^{n+1}\}\|^{2}+b)\hat{\varphi}_{\ell,m}^{n+1}-(a\|\{\hat{\psi}^{n}_{\ell,m}\}\|^{2}+b)\{\hat{\psi}^{n}_{\ell,m}\}
			-\alpha_{n}\gamma_{n}e_{0,0}=0. \quad
		\end{aligned}
	\end{equation}
	Thus, $\gamma_{n}=e_{0,0}^{\mathcal{T}}\nabla F_{h}(\{\hat{\psi}^{n}_{\ell,m}\})$.
	This yields a nonlinear equation in the form of
	\tR{
		\begin{equation}
			\begin{aligned}
				\{\hat{\varphi}^{n+1}_{\ell,m}\}&=\left( \alpha_{n}D+(a\|\{\hat{\varphi}^{n+1}_{\ell,m}\}\|^{2}+b)\mathrm{I} \right)^{-1}
				\left( (a\|\{\hat{\psi}^{n}_{\ell,m}\}\|^{2}+b)\{\hat{\psi}^{n}_{\ell,m}\}
				-\alpha_{n}\nabla F_{h}(\{\hat{\psi}^{n}_{\ell,m}\}) \right), \\
				&
				\quad\mbox{$\nabla F_{h}(\{\psi^{n}_{\ell,m}\})_{0,0}=0$.}
			\end{aligned}
			\label{eq:D_LB_model_G2_2}
	\end{equation}}
	We solve such a fixed point problem by Newton's method.
	\par
	Algorithm \ref{alg:Nesterov} and Algorithm \ref{alg:BPGK} summarize the procedures of the Nesterov, ANesterov and AA-BPG-M methods for the spherical LB model. 
	The step size estimations are detailed in Algorithm \ref{alg:line-search2},
	which requires a sufficient decrease of the energy function $E_{h}$ compared with the traditional line search in Algorithm \ref{alg:line-search}.
	It is important to note that unlike the AA-BPG algorithms, the Nesterov and ANesterov algorithms do not include a restart step for energy dissipation.
	\begin{algorithm}[!htp]
		\caption{Nesterov and ANesterov methods}
		\label{alg:Nesterov}
		\begin{algorithmic}[1]
			\REQUIRE Initial values: \;$\{\hat{\varphi}^{1}_{\ell,m}\}=\{\hat{\varphi}^{0}_{\ell,m}\}$, $\alpha>0,\; \alpha_{0}>0,\; \alpha_{\min}>0,\; w_{0}>0,\; \eta>0,\; \bar{w}>0$ 
			and $\tau>0$, $n=1$
			\WHILE{$\|\nabla E_{h}(\{\hat{\varphi}^{n}_{\ell,m}\})\|<\tau$}
			\STATE Introduce an auxiliary variable:\; $\hat{\psi}^{n}_{\ell,m}=\hat{\varphi}^{n}_{\ell,m}-w_{n}(\hat{\varphi}^{n}_{\ell,m}-\hat{\varphi}^{n-1}_{\ell,m})$
			\IF{Nesterov method} 
			\STATE Choose a fixed step size $\alpha_{n}=\alpha$;\; update $z^{n}_{\ell,m}$ by Eq.\;\eqref{eq:D_LB_model_G2_1}
			\ENDIF
			\IF{ANesterov method} 
			\STATE Calculate $\alpha_{n}$ by the Algorithm \ref{alg:line-search2};\; update $z^{n}_{\ell,m}$ by Eq.\;\eqref{eq:D_LB_model_G2_1}
			\ENDIF
			\IF{$E_{h}(\{\hat{\varphi}^{n}_{\ell,m}\})-E_{h}(\{z^{n}_{\ell,m}\})\geq  \eta\|\{\hat{\varphi}^{n}_{\ell,m}\}-\{z^{n}_{\ell,m}\}\|^{2}$ or $\alpha_{n}<\alpha_{\min}$}
			\STATE $\hat{\varphi}^{n+1}_{\ell,m}=z^{n}_{\ell,m}$,\; update $w_{n+1}\in[0,\bar{w}]$
			\ENDIF
			\STATE $n=n+1$
			\ENDWHILE
		\end{algorithmic}  
	\end{algorithm}
	
	\begin{algorithm}[!htb]
		\caption{AA-BPG-M methods}
		\label{alg:BPGK}
		\begin{algorithmic}[1]
			\REQUIRE Initial value: \;$\{\hat{\varphi}^{1}_{\ell,m}\}=\{\hat{\varphi}^{0}_{\ell,m}\}$, $\alpha_{0}>0,\; \alpha_{\min}>0,\; w_{0}>0,\; \eta>0,\; \bar{w}>0$ 
			and $\tau>0$, $n=1$
			\WHILE {$\|\nabla E_{h}(\{\hat{\varphi}^{n}_{\ell,m}\})\|<\tau$}
			\STATE Introduce an auxiliary variable:\; $\hat{\psi}^{n}_{\ell,m}=\hat{\varphi}^{n}_{\ell,m}-w_{n}(\hat{\varphi}^{n}_{\ell,m}-\hat{\varphi}^{n-1}_{\ell,m})$
			\IF{AA-BPG-2 method} 
			\STATE Calculate $\alpha_{n}$ by the Algorithm \ref{alg:line-search2};\; update $z^{n}_{\ell,m}$ by Eq.\;\eqref{eq:D_LB_model_G2_1}
			\ENDIF
			\IF{AA-BPG-4 method} 
			\STATE Calculate $\alpha_{n}$ by the Algorithm \ref{alg:line-search2};\; update $z^{n}_{\ell,m}$ by Eq.\;\eqref{eq:D_LB_model_G2_2}
			\ENDIF
			\IF{$E_{h}(\{\hat{\varphi}^{n}_{\ell,m}\})-E_{h}(\{z^{n}_{\ell,m}\})\geq  \eta\|\{\hat{\varphi}^{n}_{\ell,m}\}-\{z^{n}_{\ell,m}\}\|^{2}$ or $\alpha_{n}<\alpha_{\min}$}
			\STATE $\hat{\varphi}^{n+1}_{\ell,m}=z^{n}_{\ell,m}$,\; update $w_{n+1}\in[0,\bar{w}]$
			\ELSE
			\STATE $\hat{\varphi}^{n+1}_{\ell,m}=\hat{\varphi}^{n}_{\ell,m}$,\; set $w_{n+1}=0$
			\ENDIF
			\STATE $n=n+1$
			\ENDWHILE
		\end{algorithmic}  
	\end{algorithm}
	
	\begin{algorithm}[!htp]
		\caption{Estimate step size $\alpha_{n}$ at $\{\hat{\psi}^{n}_{\ell,m}\}$}
		\label{alg:line-search2}
		\begin{algorithmic}[1]
			\REQUIRE $\{\hat{\psi}^{n}_{\ell,m}\}$, $\eta>0$ and $\rho\in(0,1)$ and $\alpha_{\min},\; \alpha_{\max}>0$
			\STATE Initialize $\alpha_{n}$ by BB step \eqref{eq:BB_step}
			\FOR{$i=1,2,\cdots$}
			\IF{AA-BPG-2 and ANesterov} 
			\STATE Update $z^{n}_{\ell,m}$ by Eq.\;\eqref{eq:D_LB_model_G2_1}
			\ENDIF
			\IF{AA-BPG-4} 
			\STATE Update $z^{n}_{\ell,m}$ by Eq.\;\eqref{eq:D_LB_model_G2_2}
			\ENDIF \IF{$E_{h}(\{\hat{\psi}^{n}_{\ell,m}\})-E_{h}(\{z^{n}_{\ell,m}\})\geq  \eta\|\{z^{n}_{\ell,m}\}-\{\hat{\psi}^{n}_{\ell,m}\}\|^{2}$ or $\alpha_{n}<\alpha_{\min}$}
			\STATE Break
			\ELSE
			\STATE $\alpha_{n+1}=\rho\alpha_{n}$
			\ENDIF
			\ENDFOR
			\STATE Output $\alpha_{n}=\max(\min(\alpha_{n},\alpha_{\max}),\alpha_{\min})$
		\end{algorithmic}  
	\end{algorithm}
	
	\par
	Jiang et al. \cite{JiangNBPG} have demonstrated the convergence of the AA-BPG-M algorithms for the LB model in a Euclidean domain, which is independent on the spatial discretization.
	As a result, we can directly establish the convergence properties of the AA-BPG-M methods for the spherical LB model, as shown in
	Theorem \ref{theorem:convergence_result_AA-BPG}.
	This result holds for the Nesterov and ANesterov schemes,
	as long as their sequences maintain energy dissipation.
	
	\begin{theorem}
		\label{theorem:convergence_result_AA-BPG}
		Let $E_{h}(\hat{\varphi})=G_{h}(\hat{\varphi})+F_{h}(\hat{\varphi})$ be the spherical LB energy function, and the sequence $\{ \hat{\varphi}^{n}\}$ is generated by Algorithm \ref{alg:BPGK}, then we have
		\begin{enumerate}
			\item For the AA-BPG-2 method, if $\{ \hat{\varphi}^{n}\}$ is bounded, 
			then $\{ \hat{\varphi}^{n}\}$ converges to 
			some $\hat{\varphi}^{*}$ with \fR{$\nabla E_{h}(\hat{\varphi}^{*})=0$.} 
			\item For the AA-BPG-4 method, $\{ \hat{\varphi}^{n}\}$ converges to 
			some $\hat{\varphi}^{*}$ with \fR{$\nabla E_{h}(\hat{\varphi}^{*})=0$.} 
		\end{enumerate} 
	\end{theorem}

	\subsection{Principal mode analysis (PMA) method}
	\label{subsec:PMAM}
	Besides iterative methods,  the investigation of good initial values is crucial for accelerating the process of finding the desired stationary states of the non-convex and non-linear optimization problem \eqref{eq:D_LB_model}.
	In the spherical harmonic pseudo-spectrum method, estimating initial values is to give weightings of $(N+1)^{2}$ basis functions.
	Generally, it is challenging to select appropriate weightings without any prior knowledge of the desired ordered structure in such a multi-solution problem.
	Fortunately, due to the completeness of the discrete spherical harmonic expansion \eqref{eq:finite_SHT} in $L^{2}$ space, 
	the decay rate of spherical harmonic coefficients is $o(N^{-p})$ \cite{2012Spherical}, where $p$ denotes the smoothness of the stationary state.
	This implies that only a few basis functions play a dominant role in the configuration of the stationary state.
	Therefore, we define the basis functions with the first few larger amplitudes as the principal spherical harmonics, 
	and refer to the corresponding spherical degree $\ell$ and order $m$ as the principal mode numbers and principal mode directions, respectively.
	To speed up the iteration of the optimization algorithms, we propose the PMA method to estimate initial values for the spherical LB model.
	\par
	The PMA method uses principal spherical harmonics to capture the primary characteristics of the stationary states in the beginning.
	We now provide explicit formulations for selecting principal spherical harmonics based on the generic features of the model and the symmetries of the desired ordered structures.
	First, we get the principal mode number $\ell$ by analytically analyzing the essential feature of spherical LB free energy \eqref{eq:LB_model}.
	According to the spherical harmonic expansion, the spherical LB free energy  can be expressed as
	\begin{equation}
		\begin{aligned}
			E(\{\hat{\varphi}_{\ell,m}\})=\underbrace{\frac{\xi^{2}}{2}\sum_{\ell,m}\left(1-\ell(\ell+1)/R^2 \right)^{2}|\hat{\varphi}_{\ell,m}|^{2}}_{G(\{\hat{\varphi}_{\ell,m}\})}+F(\{\hat{\varphi}_{\ell,m}\}).
		\end{aligned}
		\label{eq:PMA_1}
	\end{equation}
	In the LB model, the  ordered states can occur from homogeneous state, initially near a critical wave vector \cite{1975Brazovskii}.
	An important observation is that by considering $\xi\rightarrow \infty$,
	in this case,
	the primary wave vectors $\mbk$ are restricted to the modes lying on the circles $|\mbk|=1$ and $|\mbk|=q$ to prevent free energy from growing indefinitely \cite{Stability_LP_Jiang}.
	Similarly, in our context, $G(\{\hat{\varphi}_{\ell,m}\})$ contributes more to free energy as $\xi\rightarrow \infty$ compared with the high order non-liner term $F(\{\hat{\varphi}_{\ell,m}\})$.
	Since the equilibrium state in a physical system has finite free energy, 
	we can deduce from $\sum_{\ell,m}|\hat{\varphi}_{\ell,m}|^{2}<+\infty$ that
	\begin{equation}
		\begin{aligned}
			1-\ell(\ell+1)/R^2 =0.
		\end{aligned}
		\label{eq:optimal_degree}
	\end{equation}
	Therefore, for a given sphere radius $R=\sqrt{\ell_{0}(\ell_{0}+1)}$, the single degree $\ell_{0}$ 
	is the dominant mode number that minimizes the energy value and facilitates the formation of ordered structures.
	The initial value can be written as
	\begin{equation}
		\begin{aligned}
			\varphi(\theta,\phi)=\sum^{\ell_{0}}_{m=-\ell_{0}} \hat{\varphi}_{\ell_{0},m}Y^{m}_{\ell_{0}}(\theta,\phi).
		\end{aligned}
		\label{eq:stable_phase}
	\end{equation}
	\par
	Second, we determine the principal mode directions $m$ in \eqref{eq:stable_phase} by the desired symmetry of equilibrium structures.
	To  construct a good initial value for the desired phase, we utilize the relations among the spherical symmetry, subgroups of $O(3)$ and spherical harmonics \cite{2018_Sanjay_symmetry}.
	Let $(x,y,z)\in \mathbb{R}^{3}$, $r^{2}=x^{2}+y^{2}+z^{2}$ and $\hat{x}, \hat {y}, \hat{z}$ be the directional derivatives
	satisfying $\hat{x}^{2}+\hat{y}^{2}+\hat{z}^{2}=0$. 
	With $\hat{\xi}:=\hat{x}-i\hat{y}$ and
	$\hat{\eta}:=\hat{x}+i\hat{y}$,
	these operators fulfil
	\begin{equation}
		\begin{aligned}
			& \hat{\xi}\hat{\eta}(\frac{1}{r})\lvert _{r=1}=-\hat{z}^{2},\\
			&\hat{z}^{(\ell-m)}(\hat{\xi^{m}}+\hat{\eta}^{m})(\frac{1}{r})\lvert_{r=1}=(-1)^{\ell-m}\sqrt{2(\ell-m)! (\ell+m)!} \Re {Y_{\ell}^{m}},\\
			& \hat{z}^{\ell}(\frac{1}{r})\lvert_{r=1} =(-1)^{\ell} \ell! Y_{\ell}^{0}, \\
			& i\hat{z}^{(\ell-m)}(\hat{\xi^{m}}-\hat{\eta}^{m})(\frac{1}{r})\lvert_{r=1}=(-1)^{\ell-m}\sqrt{2(\ell-m)! (\ell+m)!} \Im {Y_{\ell}^{-|m|}}.
		\end{aligned}
		\label{eq:PMA_relationship}
	\end{equation}
	Table \ref{tab:PMA_subgroup} summarizes the relationships between spherical harmonics and symmetric subgroups in $O(3)$. 
	\begin{table}[!htb]
		\centering
		\caption{Subgroups of $O(3)$ and their spherical harmonics of degree $\ell$.
			Here $s\in \{ 0,1\}$ and $p,q\in \mathbb{N}\cup \{0\}$}  
			\begin{tabular}{|p{1.6cm}|p{3.5cm}|p{2.2cm}|} 
				\hline
				subgroup &  spherical harmonics & degree $\ell$   \\
				\hline
				$\mathbb{T}$ & $\mathcal{T}^{s}_{6} \mathcal{T}^{p}_{4} \mathcal{T}_{3}^{q}(1/r)\lvert_{r=1}$ & $6s+4p+3q$ \\
				\hline
				$\mathbb{O}$ & $\mathcal{O}^{s}_{9} \mathcal{O}^{p}_{6} \mathcal{O}_{4}^{q}(1/r)\lvert_{r=1}$ & $9s+6p+4q$ \\
				\hline
				$\mathbb{I}$ & $\mathcal{I}^{s}_{15} \mathcal{I}^{p}_{10} \mathcal{I}_{6}^{q}(1/r)\lvert_{r=1}$ & $15s+10p+6q$ \\
				\hline
				$Z_{n}$  & $\hat{z}^{p}\mathcal{C}_{qn}(1/r)\lvert_{r=1}$, $\hat{z}^{p}\mathcal{S}_{qn}(1/r)\lvert_{r=1}$ & $p+qn$ \\
				\hline
			\end{tabular}
			\label{tab:PMA_subgroup}
		\end{table}
		The symmetric operators in the second column are given by
		\begin{equation}
			\begin{aligned}
				\mathcal{T}_{3} &=\frac{i}{4} \hat{z}(\hat{\xi}^{2}-\hat{\eta}^{2}), \quad
				\mathcal{T}_{4} =\frac{1}{4} \left[ 14\hat{z}^{4}+(\hat{\xi}^{4}+\hat{\eta}^{4}) \right], \quad
				\mathcal{T}_{6} =\frac{1}{32} \left[ (\hat{\xi}^{6}+
				\hat{\eta}^{6})-33\hat{z}^{4}(\hat{\xi}^{2}+\hat{\eta}^{2}) \right], \\ 
				\mathcal{O}_{4} &= 14 \hat{z}^{4}+(\hat{\xi}^{4}+\hat{\eta}^{4}), \quad
				\mathcal{O}_{6} = \hat{z}^{2}(\hat{\xi}^{4}+\hat{\eta}^{4})-2\hat{z}^{6}, \quad
				\mathcal{O}_{9} = i\left[ \hat{z}(\hat{\xi}^{8}-\hat{\eta}^{8})-34\hat{z}^{5}(\hat{\xi}^{4}-\hat{\eta}^{4}) \right], \\
				\mathcal{I}_{6} &= 11\hat{z}^{6}+\hat{z}(\hat{\xi}^{5}+\hat{\eta}^{5}), \quad
				\mathcal{I}_{10} = 494\hat{z}^{10}-228\hat{z}^{5}(\hat{\xi}^{5}+\hat{\eta}^{5})+(\hat{\xi}^{10}+\hat{\eta}^{10}), \\
				\mathcal{I}_{15} &=i\left[ -10005\hat{z}^{10}(\hat{\xi}^{5}-\hat{\eta}^{5})+522\hat{z}^{5}(\hat{\xi}^{10}-\hat{\eta}^{10})+(\hat{\xi}^{15}-\hat{\eta}^{15}) \right], \\             \mathcal{C}_{n}&=\left(\hat{\xi}^{n}+\hat{\eta}^{n} \right),\quad  \mathcal{S}_{n}=i\left(\hat{\xi}^{n}-\hat{\eta}^{n} \right).
			\end{aligned}
			\label{eq:PMA_operators}
		\end{equation}
		Therefore, we can construct a good initial value with desired symmetry using Table \ref{tab:PMA_subgroup}.
		Now, let us consider an example to illustrate the utility of this approach for constructing initial values.
		For a given $\ell_{0}=6$, 
		an ordered structure with $\bbI$ symmetry requires the third column of Table \ref{tab:PMA_subgroup} to be $6$,
		that is 
		\begin{equation*}
			15s+10p+6q=6, \quad \quad s=0,\; p=0,\; q=1. 
		\end{equation*}
		Thus the initial value is
		\begin{equation}
			\begin{aligned}
				\mathcal{I}_{6}^{1}(1/r)\lvert_{r=1}= 11\hat{z}^{6}+\hat{z}(\hat{\xi}^{5}+\hat{\eta}^{5})=11 \Re Y^{0}_{6}+ \Re Y^{5}_{6}:=span \{Y_{6}^{0}+Y^{5}_{6}\}.
			\end{aligned}
		\end{equation}
		This implies that we use $\ell_0=6$ with $m=0,5$ to construct an initial configuration and estimate $R
		= \sqrt{\ell_0(\ell_0+1)} = \sqrt{42}$ for the $\mathbb{I}$-symmetric phase.

%% file: Part_4_Numerical_results.tex
		\section{Numerical results}
		\label{sec:Analysis}
		In this section, we take spotted and striped phases as examples to demonstrate the performance of the PMA method and several optimization methods, 
		including AA-BPG-2, AA-BPG-4, Nesterov, ANesterov, AGD and ACG algorithms.
		The optimization methods are applied to calculate stationary states of the finite-dimensional spherical LB model. 
		Their efficiency is presented by comparing with the SIS and ASIS methods.
		\tR{Specifically, the SIS method updates the numerical solution of the spherical LB model by
			\begin{equation*}
				\{\hat{\varphi}^{n+1}_{\ell,m}\}=(\alpha_{n}+\mathrm{I})^{-1}(\{\hat{\varphi}^{n}_{\ell,m}\}-\alpha_{n}\nabla F_{h}(\{\hat{\varphi}^{n}_{\ell,m}\})),
			\end{equation*}
			where the step size $\alpha_{n}$ is a fixed number.
			In contrast, the ASIS method applies Algorithm \ref{alg:line-search} to adaptively update $\alpha_{n}$.}
		The step sizes $\alpha$ in Nesterov and SIS approaches are chosen to guarantee the best numerical behaviour, 
		and the step sizes $\alpha_{n}$ of others are obtained adaptively by the linear search technique.
		In the following simulations, the latitude and longitude angles are discretized by 
		$N_{\theta}=512$ Gaussian nodes and $N_{\phi}=2048$ uniform grid points, respectively.
		The maximum degree $N$ of spherical harmonics is truncated at $127$.
		For all methods, the stopping criterion is that the
		gradient error satisfies $\|\nabla E_{h} \|_{\infty} < 10^{-6}$.
		We set $\eta=1.0\times 10^{-14}$ in Algorithm \ref{alg:Nesterov} and Algorithm \ref{alg:BPGK}.
		For the line search strategies in Algorithm \ref{alg:line-search} and Algorithm \ref{alg:line-search2},
		$\rho=(\sqrt{5}-1)/2$ when line searches are less than $8$,
		otherwise, it is $0.1$.
		We use the PMA method to give initial values, and demonstrate the efficiency by comparing with random initial values.
		All codes were written in the MATLAB language without a parallel implementation.
		The SHTns(v3.5) package \cite{shtns_Ishioka_2018} is adopted to implement the discrete spherical harmonic transformation.
		Numerical experiments were performed on a workstation with a 2.40 GHz CPU 
		(i5-1135G7, 2 processors).
		
		\subsection{The efficiency of the PMA method}
		\label{subsec:rslt.pma}
		
		In this subsection, we will demonstrate the efficiency of the PMA method for estimating initial values.
		It is well-known that the stationary state is sensitive to the initial
		configuration and the sphere radius because of the existence of multiple solutions.
		To show the effectiveness of the PMA method in accelerating the process of finding the desired stationary phases, 
		we will compare it with random initial values.
		
		\subsubsection{Spotted phase}
		\label{subsubsec:rslt.pma.spot}
		
\begin{figure}[!htb]
			\centering
				\includegraphics[width=74mm]{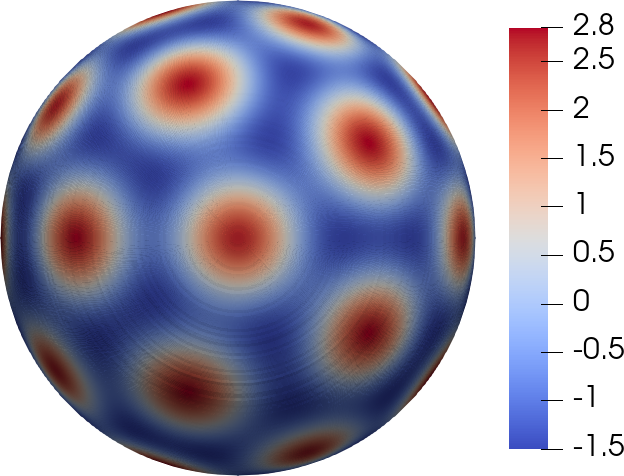}
			\caption{
				A stationary spotted phase with 32 spots when $\xi=1.0$, $\epsilon=-0.4$,
				$\lambda=0.4$.
				The initial configuration and sphere radius are given by the PMA method
			}
			\label{fig:PMA_L_10_spots_32}
		\end{figure}

		We consider a spotted phase with $32$ spots to show that the PMA method can estimate good initial values for the desired stationary spotted structures.
		Other spotted phases have the similar results.
		The model parameters are $\xi=1.0$, $\epsilon=-0.4$, $\lambda=0.4$.
		For the spotted phase, the PMA method chooses the principal mode number as $\ell=10$, and uses the $\bbI$ symmetric group to determine the principal spherical harmonics by $m=0,5,10$, i.e.,
		\begin{equation}
			\begin{aligned}
				\varphi_{S_{10}}= \hat{\varphi}_{10,0}Y^{0}_{10}
				+ \hat{\varphi}_{10,5}Y^{5}_{10}
				+ \hat{\varphi}_{10,10}Y^{10}_{10},
				\nonumber
			\end{aligned}
		\end{equation}
		where these coefficients are distributed in $(0,1]$.
		Meanwhile, the PMA method estimates sphere radius by $R=\sqrt{\ell(\ell+1)}=\sqrt{110}$.
		By choosing these initial values, a stationary $32$-spotted phase can be captured, as shown in Fig.\;\ref{fig:PMA_L_10_spots_32}.
		
		\par
		Table\;\ref{tab:PMA_L_10_spots_32} presents the success rate of different
		initial values and sphere radius in obtaining stationary states of the desired structure in $200$ experiments.
		\begin{table}[!htb]
			\caption{
				The success rate of the PMA method and random initial values to obtain the desired stationary spotted  phase with $32$ spots.
				Each case takes $200$ experiments
			}  
			\begin{center}
				\begin{tabular}{|p{3.6cm}|p{2.8cm}|p{2.0cm}|}  
					\hline  
					initial value $\varphi^{0}$   & sphere radius $R$ 
					& success rate   \\  
					\hline  
					$\varphi_{S_{10}}$		& $\sqrt{110}$		& $100\%$   \\ 
					\hline
					$\varphi_{S_{10}}$      & random number     & $4.5\%$   \\ 
					\hline
					random distribution     & $\sqrt{110}$		& $3.5\%$   \\ 
					\hline
					random distribution     & random number     & $0\%$		\\ 
					\hline  
			\end{tabular}  
		\end{center} 
		\label{tab:PMA_L_10_spots_32} 
	\end{table} 
	As the table shows, we consider four different cases.
	The initial state $\varphi_{S_{10}}$ and sphere radius $\sqrt{110}$ are given by the PMA method.
	The first row with $\varphi^0=\varphi_{S_{10}}$ and 
	$R=\sqrt{110}$ has $100\%$ success rate to obtain the desired stationary state.
	The second row uses $\varphi^0=\varphi_{S_{10}}$ but random sphere radius and
	the success rate reduces to $4.5\%$.
	If the initial value is generated randomly but the sphere radius
	is $\sqrt{110}$, the success rate becomes $3.5\%$, see the third row.
	The success rate drops to $0\%$ if both the initial configuration and the sphere radius
	are random numbers, see the fourth row.
	Therefore, the PMA method proves effective in choosing good initial states and sphere radius, significantly improving the success rate to obtain stationary states of desired spotted phases.
	
	\TR{
		\begin{remark}
			It is noted that the main work of this section shows the effect of initial values instead of the iteration algorithms.
			Here, SIS method is used to test the accuracy and efficiency of PMA for estimating initial values.
			In fact, the other optimization methods have the similar conclusion. 
			As we all known, ordered structures are deeply dependent on initial configurations and sphere radius, which means that random initial values have a sharp decreasing success rate to the desired phase.
			In contrast, we also conduct the other numerical experiments at larger $R$, which show that PMA generally has more than $50\%$ success rate much more than that of random initial values.
			The corresponding results are not presented in this paper,
			since the present numerical results are enough for demonstrating the effectiveness of the PMA for accurately estimating the initial configuration and sphere radius $R$.		
		\end{remark}	
	}

	\subsubsection{Striped phase}
	\label{subsubsec:rslt.pma.strip}
\begin{figure}[!htb]
		\centering
		\includegraphics[width=74mm]{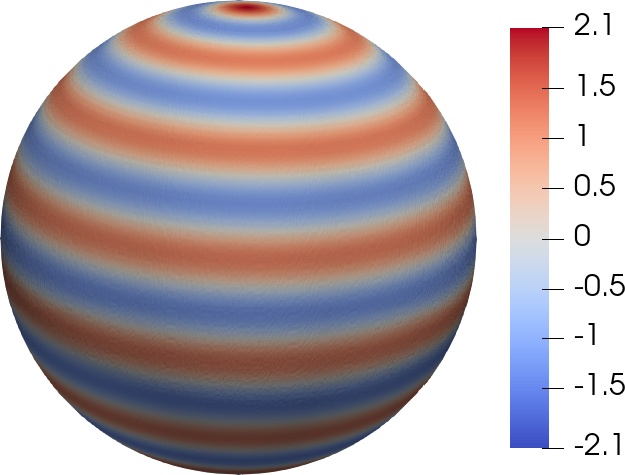}
		\caption{A stationary striped phase with \tR{$16$} strips when $\xi=1.0$,
			$\epsilon=-0.2$, $\lambda=0.0$.
			The initial configuration and sphere radius are given by the PMA method
		}
		\label{fig:PMA_L_15_strip}
	\end{figure}

	Here we take a striped phase with $16$ stripes as the desired structure to further demonstrate the efficiency of the PMA method.
	The model parameters are $\xi=1.0$, $\epsilon=-0.2$, $\lambda=0.0$.
	To obtain a $16$-striped phase with $Z_{15}$ symmetry,
	the PMA method gives the sphere radius by $R=\sqrt{240}$ and the initial state by
	\begin{equation}
		\begin{aligned}
			\varphi_{L_{15}} = \hat{\varphi}_{15,0}Y^{0}_{15}.
			\nonumber
		\end{aligned}
	\end{equation}
	Such an initial state can converge to the desired stationary $16$-striped phase, as shown in Fig.\;\ref{fig:PMA_L_15_strip}.
	This striped phase has $8$ red circles and $8$ blue circles.
	\par
	Table\;\ref{tab:PMA_L_15_strip} compares the success rate of the PMA method and random initial values in finding the desired striped phase.
	As shown in the first row, the success rate is $100\%$ when the initial values and sphere radius given by the PMA method.
	When $R=\sqrt{240}$ but the initial state is a random distribution, the success
	rate drops rapidly to $2\%$, see the third row.
	However, if the sphere radius is a random number, the success rate becomes $0\%$, 
	as shown in the second and fourth rows.
	These results show that the PMA method is also an efficient method to estimate
	initial values and sphere radius for striped phases.
	\begin{table} [!htb]
		\caption{
			The success rate of the PMA method and random initial values to obtain the desired
			structure that is a striped phase with $16$ strips.
			Each case takes $200$ experiments			
		}  
		\begin{center}
			\begin{tabular}{|p{3.6cm}|p{2.8cm}|p{2.0cm}|}  
				\hline  
				initial value $\varphi^{0}$   & sphere radius $R$ 
				&success rate   \\  
				\hline  
				$\varphi_{L_{15}}$		&$\sqrt{240}$  & $100\%$    \\ 
				\hline
				$\varphi_{L_{15}}$		&random number        & $0\%$     \\ 
				\hline
				random distribution     &$\sqrt{240}$  & $2\%$    \\ 
				\hline
				random distribution     &random  number      & $0\%$  \\    
				\hline  
		\end{tabular}  
	\end{center} 
	\label{tab:PMA_L_15_strip} 
\end{table}

\subsection{The efficiency of optimization methods}
\label{subsec:rslt.opt}
In this subsection, we show the performance of the developed optimization algorithms in computing spotted and striped phases, 
which is measured by the iterations and CPU time required to obtain the equilibrium states.
Since Sect. \ref{subsec:rslt.pma} has shown that the PMA method can estimate good initial states and sphere radius to accelerate the iterative process, 
all simulations performed below use initial values given by PMA.

\subsubsection{Spotted phase}
\label{subsub:rslt.opt.spot}

\begin{figure}[!htb]
	\centering
	\subfigure[Initial configuration $\varphi_{S_{15}}$]{
		\includegraphics[width=70mm]{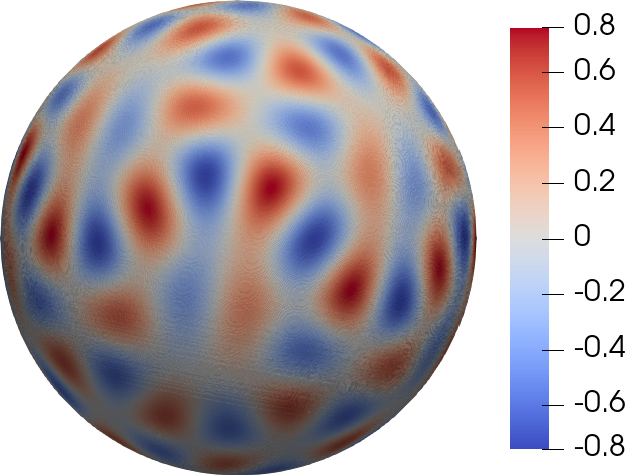}}
	\hfill
	\subfigure[Stationary spotted phase]{
		\includegraphics[width=70mm]{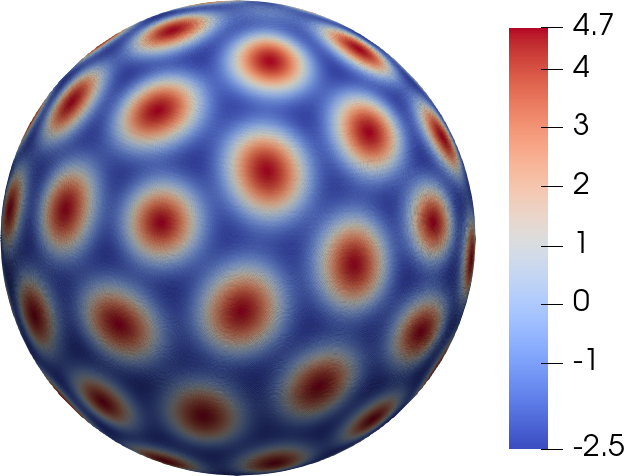}}
	\caption{
		Initial and stationary structures for the spotted phase with $60$ spots.
		$\xi=1.0$, $\epsilon=-1.0$, $\lambda=0.8$
	}
	\label{fig:AA_BPG_stationary_15_spots}
\end{figure}

First we use spotted phases to demonstrate the performance of the optimization algorithms.
By applying the PMA method, the sphere radius is $R=\sqrt{240}$ and the initial state
$\varphi^{0}$ is
\begin{equation}
	\begin{aligned}
		\varphi_{S_{15}}(\bm{r})
		=\hat{\varphi}_{15,-15}Y^{-15}_{15}
		+\hat{\varphi}_{15,-10}Y^{-10}_{15}+ \hat{\varphi}_{15,-5}Y^{-5}_{-15}
		,\nonumber
	\end{aligned}
\end{equation}
where the amplitudes are non-zero numbers and satisfy $\|\varphi^{0}\|=1$. 
Figure\;\ref{fig:AA_BPG_stationary_15_spots} presents the initial and
stationary states of a spotted phase with $60$ spots.
\par
We apply all optimization methods mentioned in Sect. \ref{sec:Numer}, as well as the SIS and ASIS methods, to compute such a spotted phase.
For fair comparison, we select the same initial state $\varphi_{S_{15}}$ and sphere radius $R=\sqrt{240}$ for all algorithms.
Furthermore, all parameters in these approaches are chosen to achieve the best performance.
Concretely, in AA-BPG-4 method \eqref{eq:convex_function_c}, we set $a=0.01$ and $b=1.0$.
For the SIS method, We choose $\alpha=0.6$, while $\alpha=0.8$ for the Nesterov method.
The AGD and ACG methods have $\alpha_{0}=0.002$,
$\alpha_{\min}=1.0\times 10^{-5}$ and $\alpha_{\max}=5.0$.
Meanwhile, parameters in the AA-BPG algorithms are
$\alpha_{0}=0.02$, $\alpha_{\max}=5.0$ and $\alpha_{\min}=0.01$, but
$\alpha_{\max}=20.0$ in the ASIS method and $\alpha_{\min}=1.0\times 10^{-5}$ in the ANesterov method.

\par
Table\;\ref{tab:AA_BPG_S_15} shows corresponding convergent results.
\begin{table} [!htb] 
	\caption{
		Convergent results of all algorithms for computing the spotted phase with $60$
		spots. 
		The initial state is $\varphi_{S_{15}}$ with $R = \sqrt{240}$, $\xi=1.0$,
		$\epsilon=-1.0$, $\lambda=0.8$
	}  
	\begin{center}
		\begin{tabular}{|p{2.0cm}|p{1.5cm}|p{2.2cm}|p{2.0cm}|p{3.8cm}|}  
			\hline  
			Method &Iteration & CPU time (s) & $ \|\nabla E_{h}(\hat{\varphi})\|_{\infty} $
			& Equilibrium energy $E_{s}$   \\  
			\hline  
			AA-BPG-2 &  172 & 48.65 & $0.98 \times 10^{-6}$ & -4.2399690344\\ 
			\hline 
			AA-BPG-4 &  130 & 37.13 & $0.67 \times 10^{-6}$ &-4.2399690344\\ 
			\hline 		
			SIS       &  994 & 175.77& $0.99 \times 10^{-6}$ & -4.2399690344 \\ 
			\hline	
			ASIS      &  380 & 95.92 & $0.99 \times 10^{-6}$ & -4.2399690344 \\ 
			\hline 
			Nesterov & 159    & 26.93  & $0.98 \times 10^{-6}$ & -4.2399690344 \\
			\hline
			ANesterov &  605 &216.11 & $0.99 \times 10^{-6}$  &-4.2399690344 \\ 
			\hline
			ACG & 119450  &96238.07  & $0.97 \times 10^{-6}$    & -4.2254676259 \\
			\hline	
			AGD & 142050 & 36190.21 & $0.99\times 10^{-6}$ & -4.2399690344 \\
			\hline	
		\end{tabular}  
	\end{center} 
	\label{tab:AA_BPG_S_15} 
\end{table} 
From the last column of the table,
we observe that all algorithms, except for the ACG algorithm with $E_{s}=-4.2254676259$, converge to a constant energy value
$E_{s}=-4.2399690344$ indicating the same stationary structure.
The slight energy difference can be attributed to the non-linearity of spherical LB free energy, which may lead to inaccurate estimations of the iterative directions in the ACG Algorithm, ultimately preventing the free energy function from reaching the minimum value.
Furthermore, we observe that AA-BPG-2/4 and Nesterov algorithms require less than $50$ seconds of CPU time and fewer than $180$ iterations, outperforming than others.
Specifically, for such a stationary spotted phase,
the AA-BPG-4 method achieves convergence in $130$ iterations and $37.13$ seconds, slightly better than the AA-BPG-2 method.
The AA-BPG-4 method has less CPU time, and is about $4$ times faster than the SIS method, $2$ times faster than the ASIS method, $6$ times faster than the ANesterov method, 
$2500$ times faster than the ACG method and $950$ times than the AGD method.
It is worth noting that the Nesterov method costs less CPU time 
than the AA-BPG-4 method, despite requiring more iterations.
This is because the AA-BPG-2/4 methods spend additional time on the line search in each iteration.
This can also explain why the ACG method costs more CPU time than the AGD method.

\par
\begin{figure}[!htb]
	\centering
	\subfigure[ Energy difference over iterations ]{
		\includegraphics[width=74mm,height=60mm]{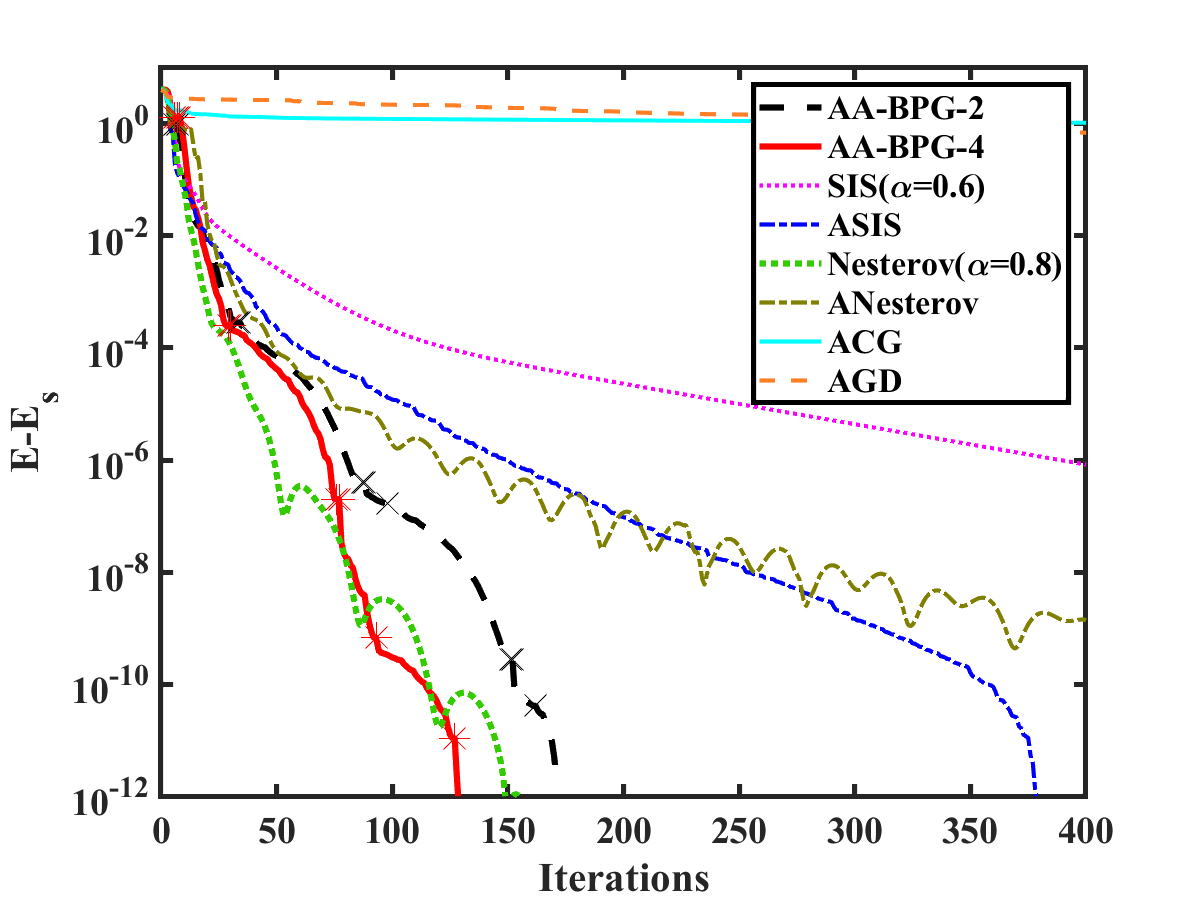}}
	\hspace{0.5cm}
	\subfigure[Energy difference over CPU time (s)]{
		\includegraphics[width=74mm,height=60mm]{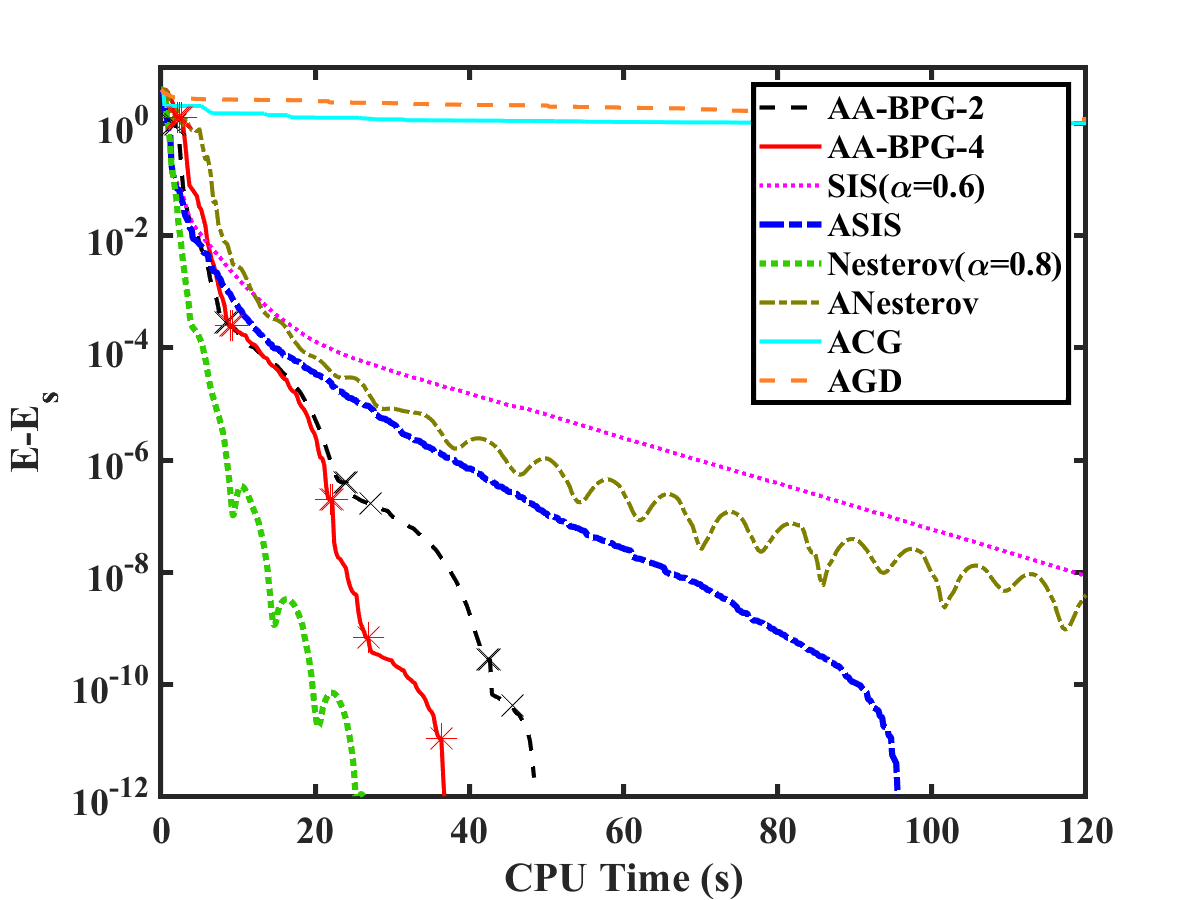}}
	\subfigure[Gradient over iterations ]{
		\includegraphics[width=74mm,height=60mm]{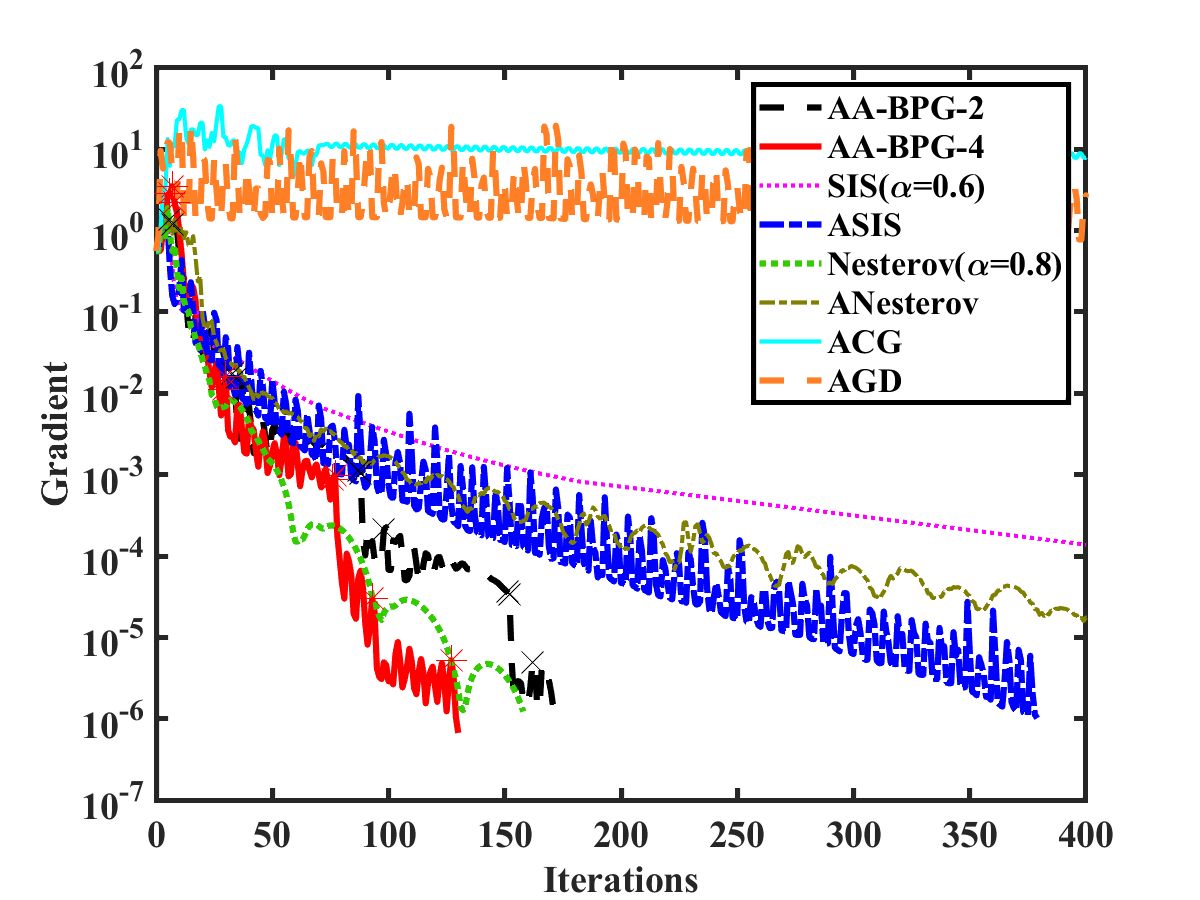}}
	\caption{Comparison of numerical behavior of the AA-BPG-2, AA-BPG-4, Nesterov,
		ANesterov, ACG and AGD methods as well as SIS, ASIS methods for computing the
		spotted phase on a sphere of radius $R=\sqrt{240}$.
		The $\times$ markers denote the restart steps in AA-BPG-2 and
		AA-BPG-4 algorithms			
	}
	\label{fig:BPGk_spot_15}
\end{figure}

Figure\;\ref{fig:BPGk_spot_15} depicts the iterative process, including  relative energy difference and the gradient error over iterations and CPU time. 
The reference energy value of the ACG algorithm is set as $E_{s}= -4.2254676259$, and the other algorithms have $E_{s}=-4.2399690344$.
From these profiles, it is evident that the AA-BPG-2, AA-BPG-4 and Nesterov methods have significantly faster convergence to equilibrium states compared with the other methods.
Moreover, the AA-BPG-2/4 methods exhibit efficient energy dissipation.
However, the Nesterov method displays sharp energy decrease but with some energy oscillations.
Similarly, the ANesterov and ACG methods, despite employing a linear search strategy, also exhibit energy oscillations.
It should be noted that the convergence properties of these algorithms are theoretically based on energy dissipation.
The Nesterov method, while computationally efficient, lacks a theory to guarantee convergence due to the presence of energy oscillations.
The same conclusion is drawn for the ANesterov and ACG algorithms.
From the above results and analysis,
it can be concluded that the AA-BPG-4 algorithm has good performance in both numerical behavior and theoretical convergence and can be regarded as the most efficient method.
\par
\begin{figure}[!htb]
	\centering
	\includegraphics[width=74mm,height=60mm]{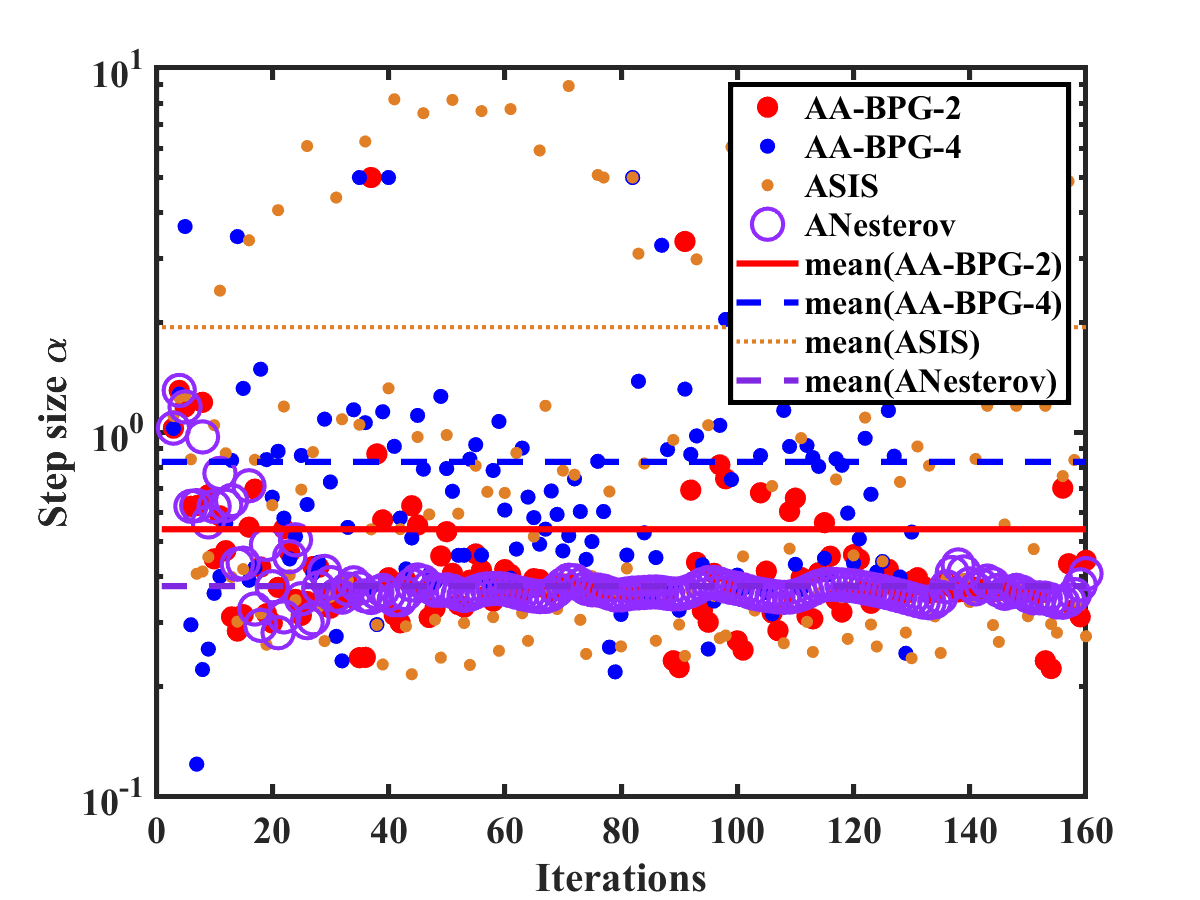}
	\hspace{0.5cm}	
	\includegraphics[width=74mm,height=60mm]{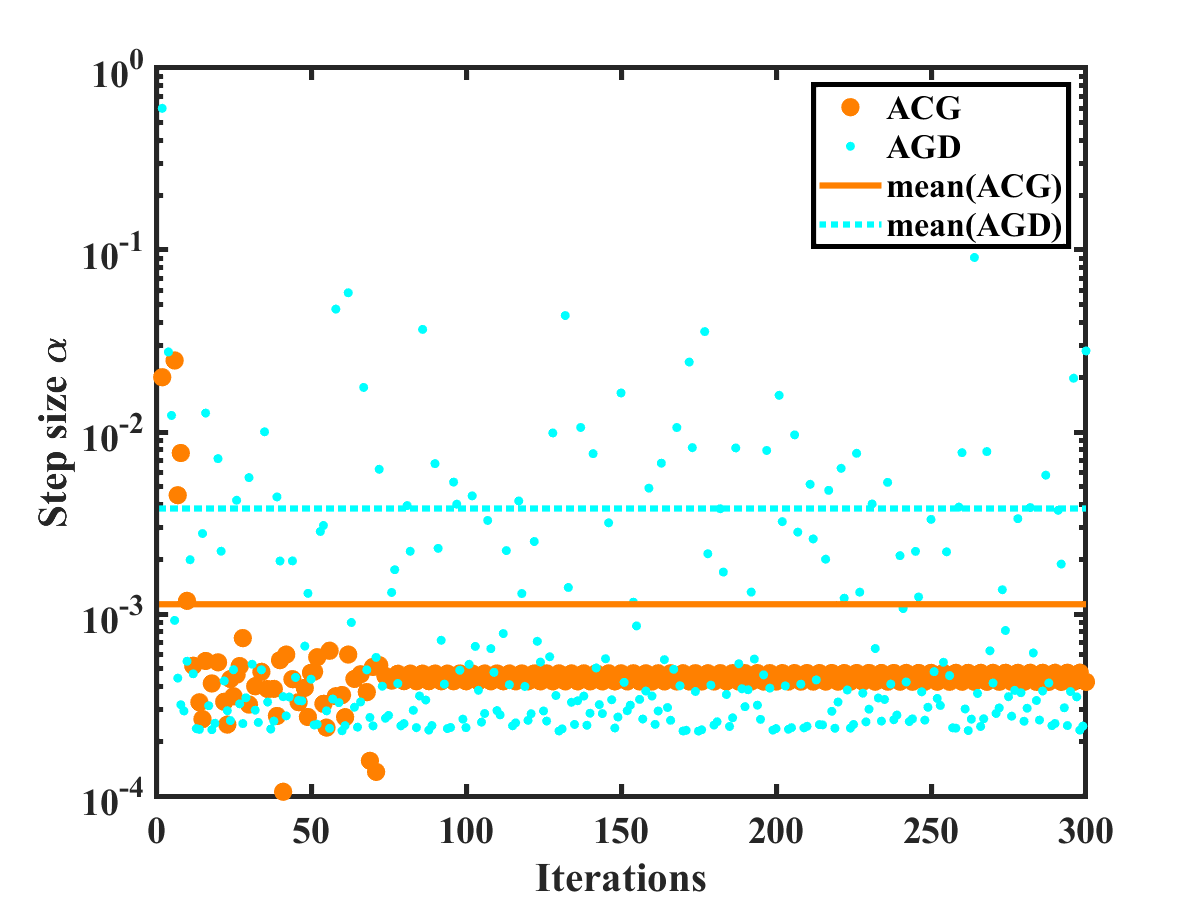}
	\caption{
		Step behavior of AA-BPG-2, AA-BPG-4 and other adaptive methods for computing the stationary spotted phase.
		The step sizes of the AGD and ACG method are $\alpha_{0}=0.002$,
		$\alpha_{\min}=1.0\times 10^{-5}$ and $\alpha_{\max}=5.0$.
		Meanwhile, the AA-BPG-2/4 methods have
		$\alpha_{0}=0.02$, $\alpha_{\max}=5.0$ and $\alpha_{\min}=0.01$, while $\alpha_{\max}=20.0$ for the ASIS method and $\alpha_{\min}=1.0\times 10^{-5}$ for the ANesterov method
	}
	\label{fig:spot_15_time_step}
\end{figure}

Figure\;\ref{fig:spot_15_time_step} presents the step sizes of adaptive schemes.
As depicted in the left figure, most of the step size in each step iteration ranges from $0.1$ to $2$ for these methods.
Specifically, the mean step sizes of AA-BPG2, AA-BPG-4, ASIS and ANesterov methods are $0.54063$, $1.5339$, $1.9404$ and $0.3772$, respectively.
In contrast, the right figure shows that the ACG and AGD method have significantly smaller step sizes.
Most of their step sizes are below $0.1$, as observed from the right figure.
The mean step sizes of the ACG and AGD methods are $0.0011$ and $0.0038$, respectively.
Importantly, these numerical results, that the ACG and AGD methods require smaller step sizes than other methods, are consistent with theoretical analysis.
\par
\begin{figure}[!htb]
	\centering
	\includegraphics[width=74mm,height=60mm]{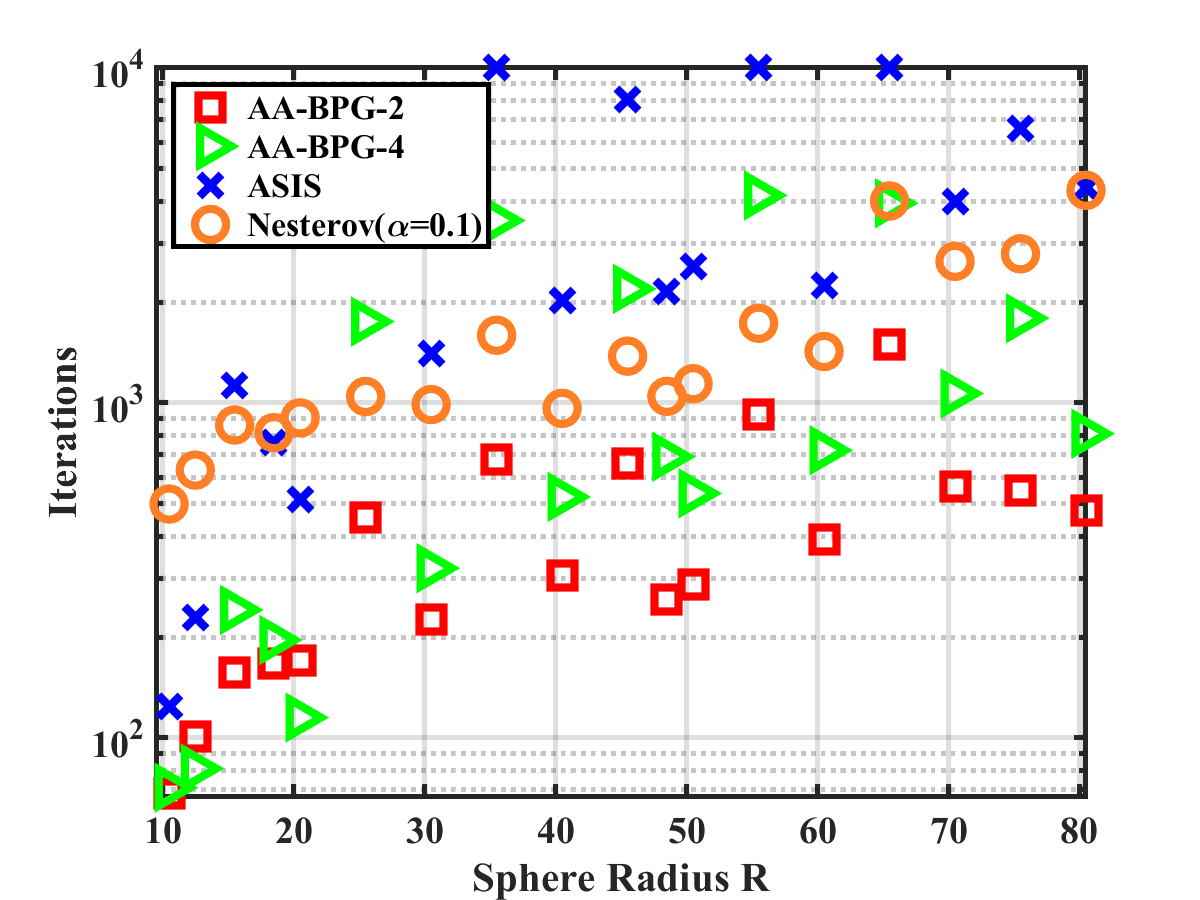}
	\hspace{0.5cm}	
	\includegraphics[width=74mm,height=60mm]{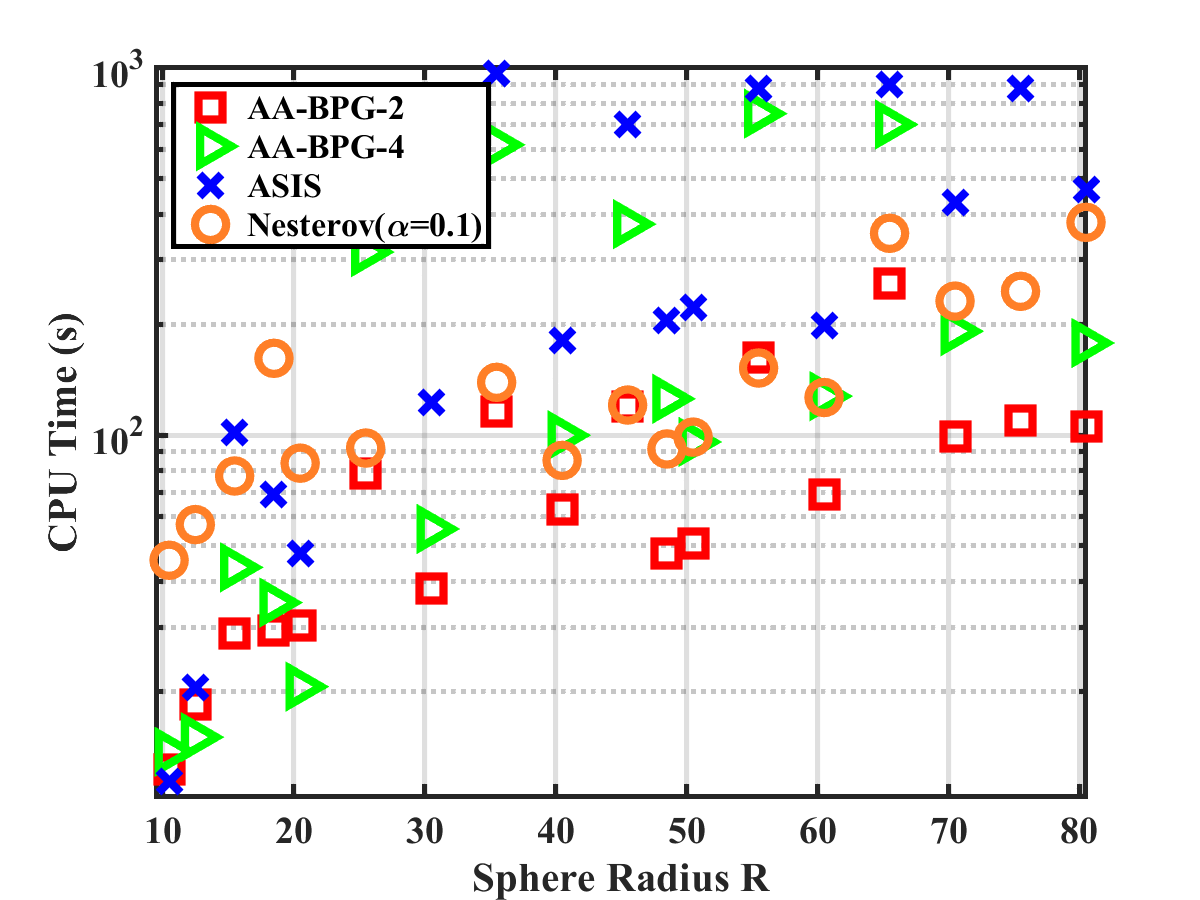}
	\caption{Required iterations and CPU time of the AA-BPG-2, AA-BPG-4, Nesterov and ASIS methods for the stationary spotted phase with different sphere radius $R$ when $\xi=1.0$, $\epsilon=-1.0$, $\lambda=0.8$.
		Here $\alpha_{0}=0.02$, $\alpha_{\min}=1.0\times 10^{-4}$ and
		$\alpha_{\max}=1.0$. 
		The initial values and sphere radius $R$ are given by
		the PMA method.
		The stopping criterion is that the gradient error satisfies $\|\nabla E_{h}
		\|_{\infty} < 10^{-6}$ or the number of iterations is greater than $1.0\times
		10^{4}$
	}
	\label{fig:ABPG_spot}
\end{figure}  
We further investigate the effectiveness of these methods with different sphere radius $R$ since the spots of stationary spotted phases increase with sphere radius $R$.
In these simulations, we fix $\alpha=0.1$ for the Nesterov method to ensure convergence.
Figure\;\ref{fig:ABPG_spot} compares the number of iterations and CPU time of the AA-BPG-2, AA-BPG-4, ASIS and Nesterov methods.
As depicted in the figures, when $R\leq 20$, the AA-BPG-2/4 methods can greatly reduce iterations and less CPU time.
However, when the sphere radius $R$ increases, the AA-BPG-2 method consistently demonstrates the least CPU time and the fewest iterations among all methods.
Based on these profiles, we conclude that the AA-BPG-2 method has the best performance in computing spotted phases.

\par
\TR{The above numerical results show that the proposed AA-BPG and Nesterov methods keep a faster convergent speed for computing the spotted phase on different spherical surfaces.
	We also compute the stationary structures at varying model parameters $\epsilon$ and $\lambda$ when $\xi=1$.
	Table \ref{tab:AA_BPG_C_15_para} compares the equilibrium energy, iterations and  CPU time of AA-BPG-2, AA-BPG-4 and ASIS methods.
	The initial conditions are set as $\alpha_{\min}=0.01$, $\alpha_{\max}=5.0$ and $\alpha_{0}=0.02$.
	The factors in AA-BPG-4 algorithms are chosen as $a=0.001$ and $b=1$.
	Obviously, compared to ASIS, AA-BPG methods can greatly reduce the iterations and CPU time in finding  stationary spotted structures in a wild range of choice of parameters.
	The performance of AA-BPG-2 and AA-BPG-4 demonstrates the certain robustness of our proposed algorithms at different model parameters.
}

\begin{table} [!htb] 
	\TR{
		\caption{
			Comparisons of numerical behevior of AA-BPG-2, AA-BPG-4 and ASIS for computing a spotted phase at varying model parameters $\epsilon$, $\lambda$.
		}  
		\begin{center}
			\begin{tabular}{|p{0.8cm}|p{0.8cm}|p{1.2cm}|p{2.0cm}|p{1.8cm}|p{2.2cm}|p{3.8cm}|}  
				\hline  
				$\epsilon$ & $\lambda$ & R	&Method &Iteration & CPU time (s) & Equilibrium energy $E_{s}$   \\  
				\hline  
				\multirow{3}*{-1} & \multirow{3}*{1} & \multirow{3}*{$\sqrt{240}$} & AA-BPG-2 &150 & 30.54 & -5.0930540417\\ 
				~ &~ & ~ & AA-BPG-4 &192 & 44.57 & -5.0930540417\\
				~&~ & ~ & ASIS & 394 & 63.81 & -5.0930540417\\
				\hline
				\hline
				\multirow{3}*{-0.75} & \multirow{3}*{1} & \multirow{3}*{$\sqrt{240}$} & AA-BPG-2 &142 & 27.57 & -3.2636101949\\
				~ &~ & ~ & AA-BPG-4 &173 & 41.18 &  -3.2636101949\\
				~ &~ & ~& ASIS &331 & 59.53 & -3.2636101949\\
				\hline
				\hline
				\multirow{3}*{-0.8} & \multirow{3}*{0.9} & \multirow{3}*{$\sqrt{240}$} & AA-BPG-2 &123 & 31.41 & -3.2357895061\\
				~ &~ & ~ & AA-BPG-4 &167 & 35.14 & -3.2357895061\\
				~ &~ & ~ & ASIS &321 & 95.38 & -3.2357895061\\
				\hline
				\hline
				\multirow{3}*{-0.3} & \multirow{3}*{0.8} & \multirow{3}*{$\sqrt{930}$} & AA-BPG-2 &211 & 44.93 & -2.2636412022\\
				~ &~ & ~ & AA-BPG-4 &226 & 53.56 & -2.2636412022\\
				~ &~ & ~ & ASIS &661 & 95.38 & -2.2636412022\\
				\hline
				\hline
				\multirow{3}*{-1} & \multirow{3}*{0.9} & \multirow{3}*{$\sqrt{930}$} & AA-BPG-2 &269 & 62.36 & -4.6793498644\\
				~ &~ & ~  & AA-BPG-4 &264 & 57.77 & -4.6793498644\\
				~ &~ & ~  & ASIS &536 &81.73 & -4.6786638033\\
				\hline
				\hline
				\multirow{3}*{-0.8} & \multirow{3}*{1} & \multirow{3}*{$\sqrt{930}$} & AA-BPG-2 & 208& 44.59 & -3.6424944647\\
				~ &~ & ~  & AA-BPG-4 &191 & 42.12  & -3.6442076494\\
				~ &~ & ~  & ASIS &370 &57.65 & -3.6424944647\\
				\hline
			\end{tabular}
		\end{center} 
		\label{tab:AA_BPG_C_15_para} 
	}
\end{table}

\subsubsection{Striped phase}
\label{subsub:rslt.opt.strip}
\begin{figure}[!htb]
	\centering
	\subfigure[Initial configuration $\varphi_{L_{60}}$]{
		\includegraphics[width=70mm]{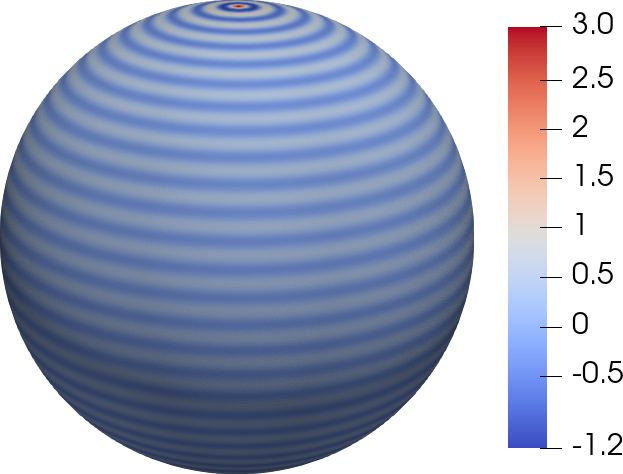}}
	\hfill
	\subfigure[Stationary striped phase]{
		\includegraphics[width=70mm]{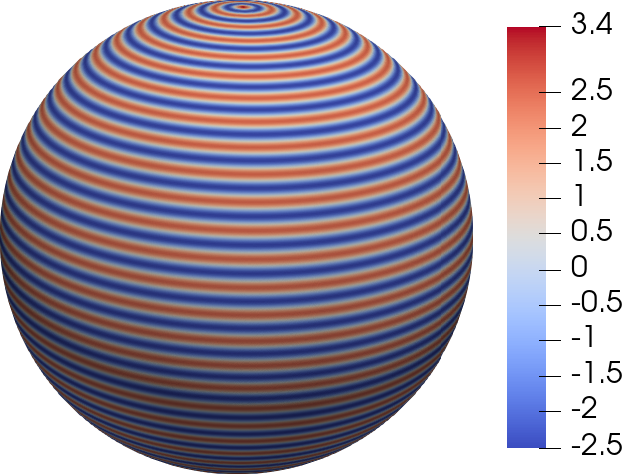}}
	\caption{
		Initial and stationary structures for the striped phase with $61$ stripes. $\xi=1.0$, $\epsilon=-0.8$, $\lambda=0.0$
	}
	\label{fig:AA_BPG_stationary_61strips}
\end{figure}

In this subsection, we focus on the efficiency of optimization methods for computing striped phases.
The model parameters are $\xi=1$, $\epsilon=-0.8$, $\lambda=0.0$.
We take a stationary striped phase with $61$ stripes as an example.
To obtain the desired striped phase, the PMA method gives $R=\sqrt{3660}$ and chooses the initial value $\varphi^{0}$ as 
\begin{equation}
	\begin{aligned}
		\varphi_{L_{60}}(\bm{r})=Y^{0}_{60}. \nonumber
	\end{aligned}
\end{equation}
Figure\;\ref{fig:AA_BPG_stationary_61strips} shows the initial and equilibrium states of the $61$-striped phase.
\par
Similarly, the above-mentioned algorithms are developed to compute a striped phase.
Here, the factors of the AA-BPG-4 method are chosen as $a=0.001$ and $b=1.0$. 
The best-performing step sizes of the SIS and Nesterov methods are $\alpha=0.8$ and $\alpha=1.5$, respectively.
Meanwhile, we choose $\alpha_{0}=0.5$, $\alpha_{\max}=45.0$ and $\alpha_{\min}=0.01$ for AA-BPG algorithms, but $\alpha_{\min}=1.0\times 10^{-8}$ for ANesterov algorithm.
The parameters in ASIS algorithm are
$\alpha_{0}=0.8$, $\alpha_{\min}=0.2$ and $\alpha_{\max}=350.0$,
while ACG algorithm has $\alpha_{0}=0.8$, $\alpha_{\min}=1.0\times 10^{-8}$ and $\alpha_{\max}=10.0$ but AGD algorithm has $\alpha_{\max}=45.0$.

\par
\begin{figure}[!htb]
	\centering
	\includegraphics[width=74mm,height=60mm]{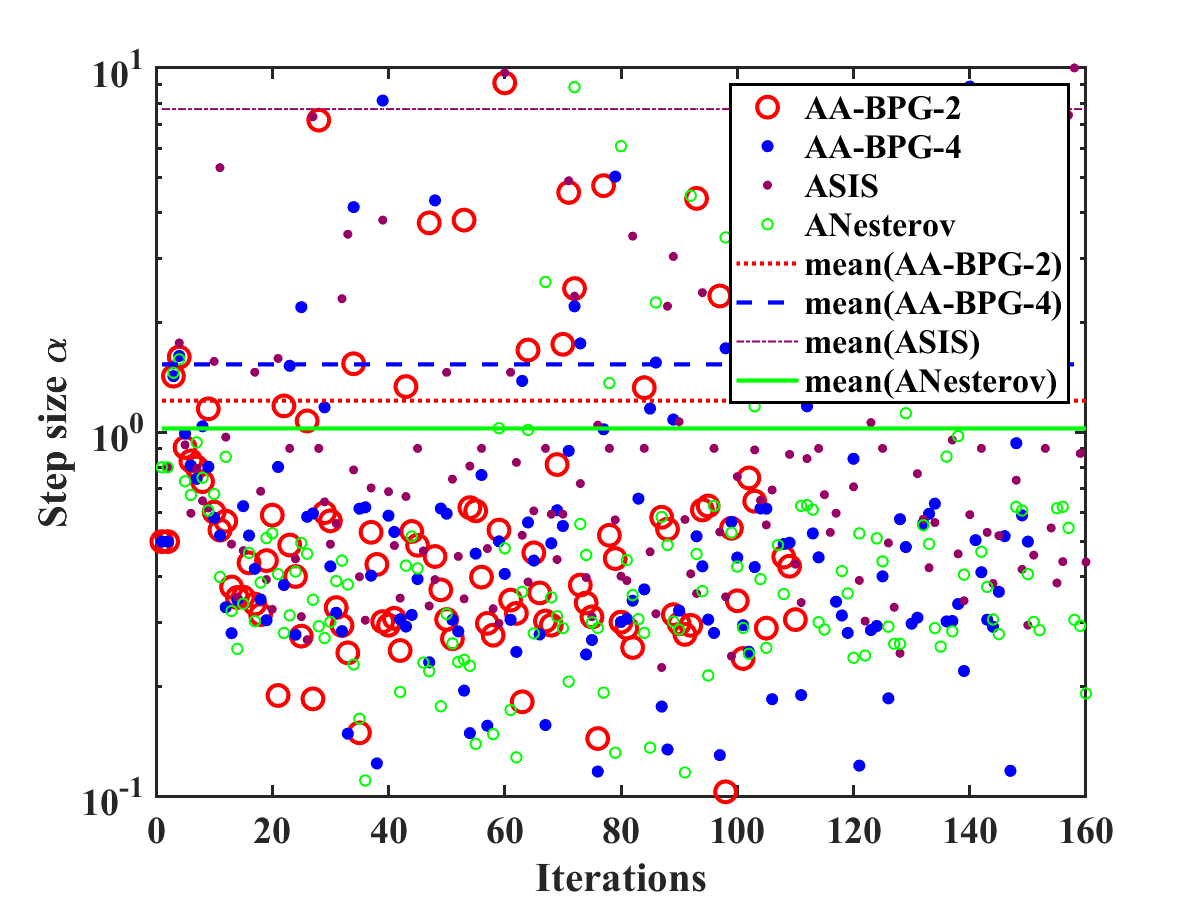}
	\hspace{0.5cm}	
	\includegraphics[width=74mm,height=60mm]{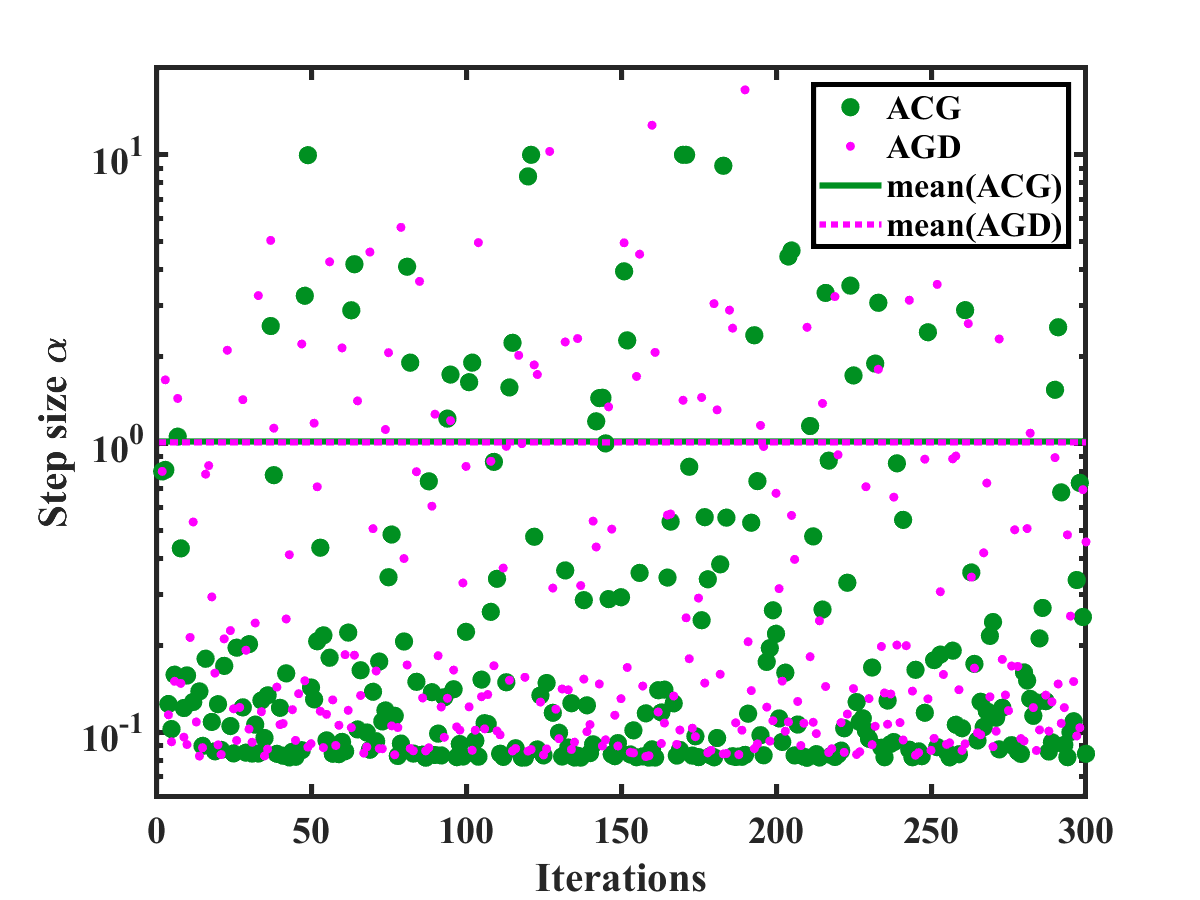}
	\caption{Step behavior of the AA-BPG-2, AA-BPG-4 methods and
		other adaptive for computing the striped phase
		Here we set $\alpha_{0}=0.5$, $\alpha_{\max}=45.0$ and $\alpha_{\min}=0.01$ for AA-BPG algorithms, but $\alpha_{\min}=1.0\times 10^{-8}$ for ANesterov algorithm.
		The ASIS algorithm has $\alpha_{0}=0.8$, $\alpha_{\min}=0.2$ and $\alpha_{\max}=350.0$,
		while ACG algorithm has $\alpha_{0}=0.8$, $\alpha_{\min}=1.0\times 10^{-8}$ and $\alpha_{\max}=10.0$ but AGD algorithm has $\alpha_{\max}=45.0$.
	}
	\label{fig:L60_time_step}
\end{figure}

Table\;\ref{tab:AA_BPG_L_60} compares convergent results of all methods.
\begin{table} [!htb] 
	\caption{Convergent results of all algorithms for computing the striped phase with $61$ stripes.
		Here $R=\sqrt{3660}$, $\xi=1.0$, $\epsilon=-0.8$, $\lambda=0.0$ and the initial state is $\varphi_{L_{60}}$
	}  
	\begin{center}
		\begin{tabular}{|p{2.0cm}|p{1.5cm}|p{2.2cm}|p{2.0cm}|p{3.8cm}|} 
			\hline  
			Method &Iteration & CPU time (s) & $\|\nabla E_{h}(\hat{\varphi})\|_{\infty}$ &
			Equilibrium energy $E_{s}$   \\  
			\hline  
			AA-BPG-2  &  111 & 38.02 & $0.95\times 10^{-6}$ & -2.2629509226 \\ 
			\hline 
			AA-BPG-4  &  153 & 47.99 & $0.98\times 10^{-6}$ & -2.2629509226 \\ 
			\hline 		
			SIS       &  2270& 367.23& $0.99\times 10^{-6}$ & -2.2629509226 \\ 
			\hline	
			ASIS      &  302 & 84.46 & $0.95\times 10^{-6}$ & -2.2629509226  \\ 
			\hline 
			Nesterov  &  2612&582.10 & $0.88\times 10^{-6}$ & -2.2629509227 \\ 
			\hline
			ANesterov &  599 &200.77 & $0.82\times 10^{-6}$ & -2.2629509226  \\ 
			\hline
			ACG       &  1445&577.28 & $0.99\times 10^{-6}$ & -2.2629509226 \\ 
			\hline	
			AGD       &  1765&723.71 & $0.99\times 10^{-6}$ & -2.2629509226 \\ 
			\hline	
		\end{tabular}  
	\end{center} 
	\label{tab:AA_BPG_L_60} 
\end{table} 
The constant equilibrium energy $E_{s}=-2.2629509226$ indicates that all algorithms accurately
converge to the same stationary striped phase.
\begin{figure}[!htb]
	\centering
	\subfigure[Energy difference over iterations]{
		\includegraphics[width=74mm,height=60mm]{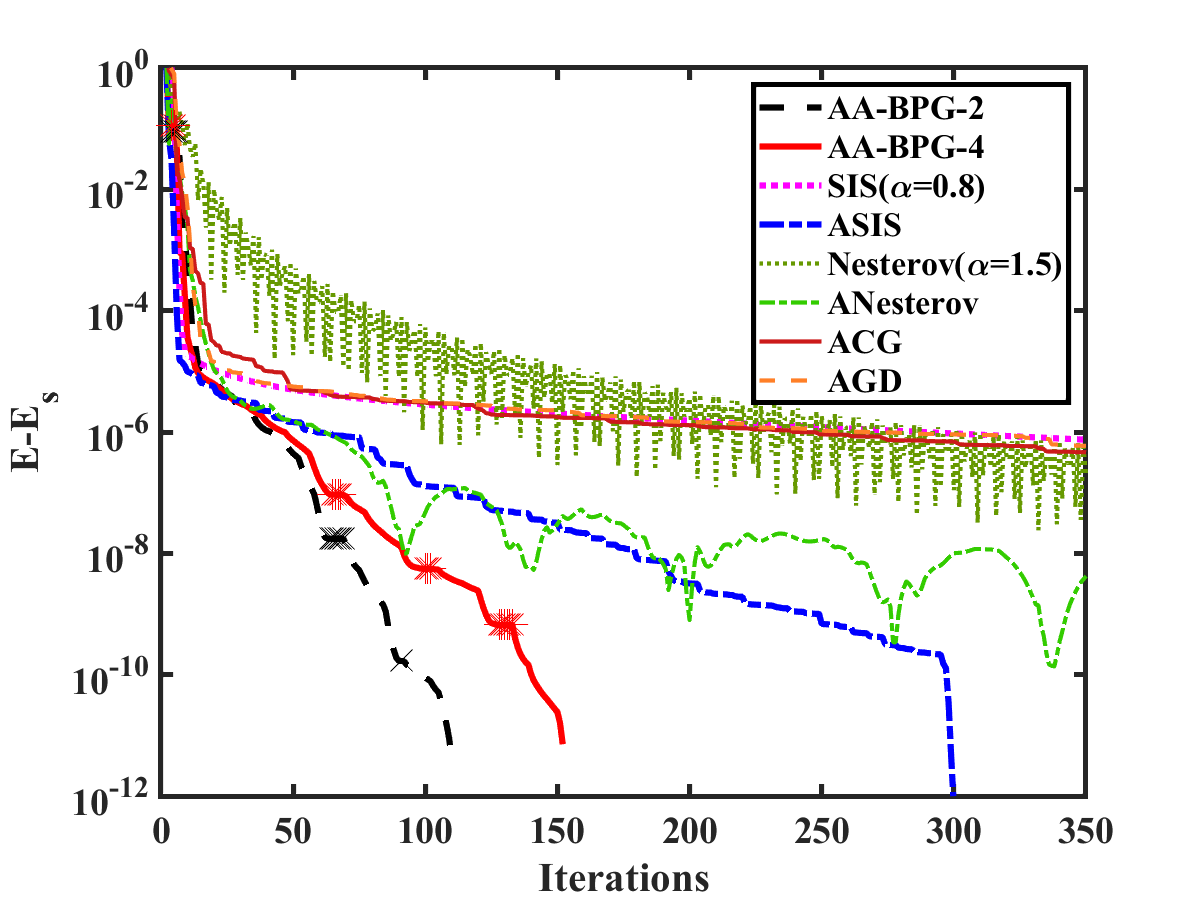}}
	\hspace{0.5cm}
	\subfigure[Energy difference over CPU time]{
		\includegraphics[width=74mm,height=60mm]{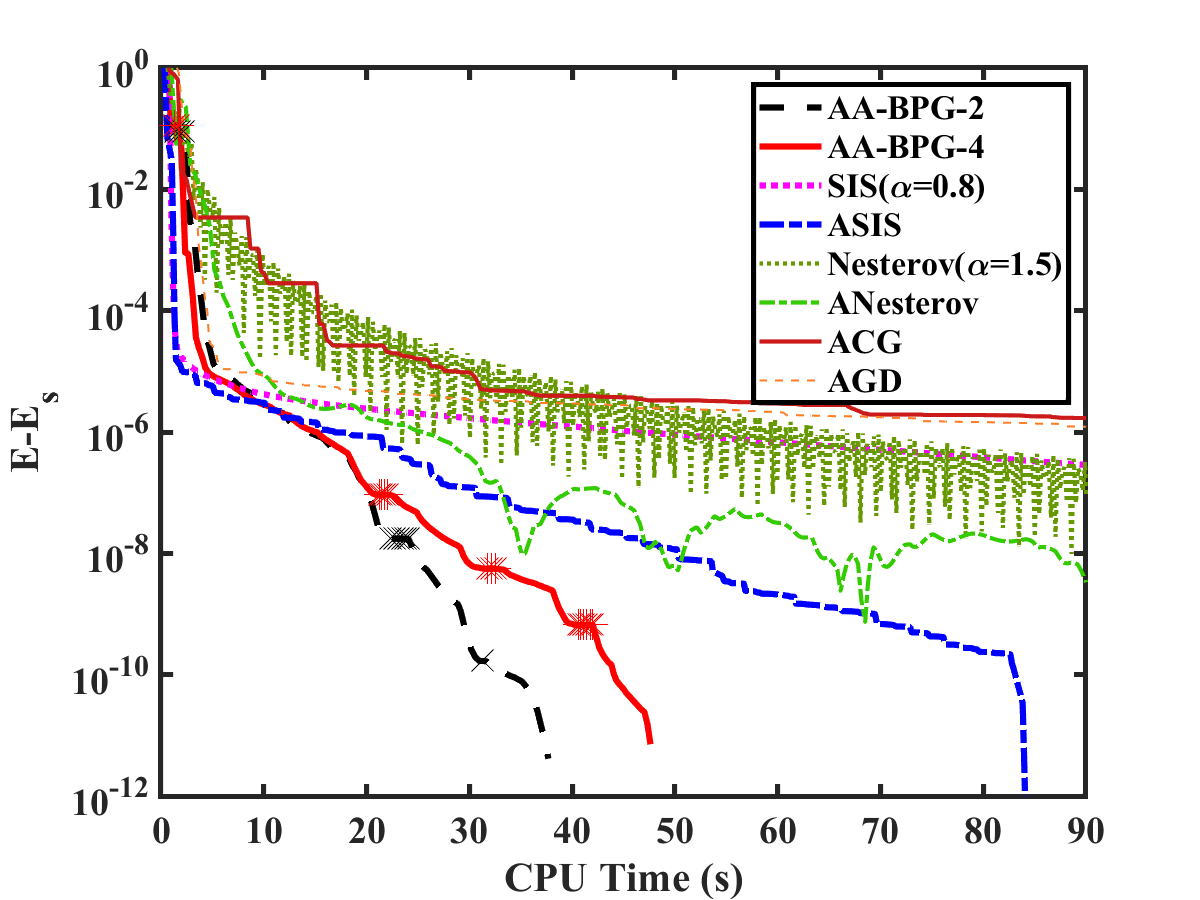}}
	\subfigure[Gradient over iterations]{
		\includegraphics[width=74mm,height=60mm]{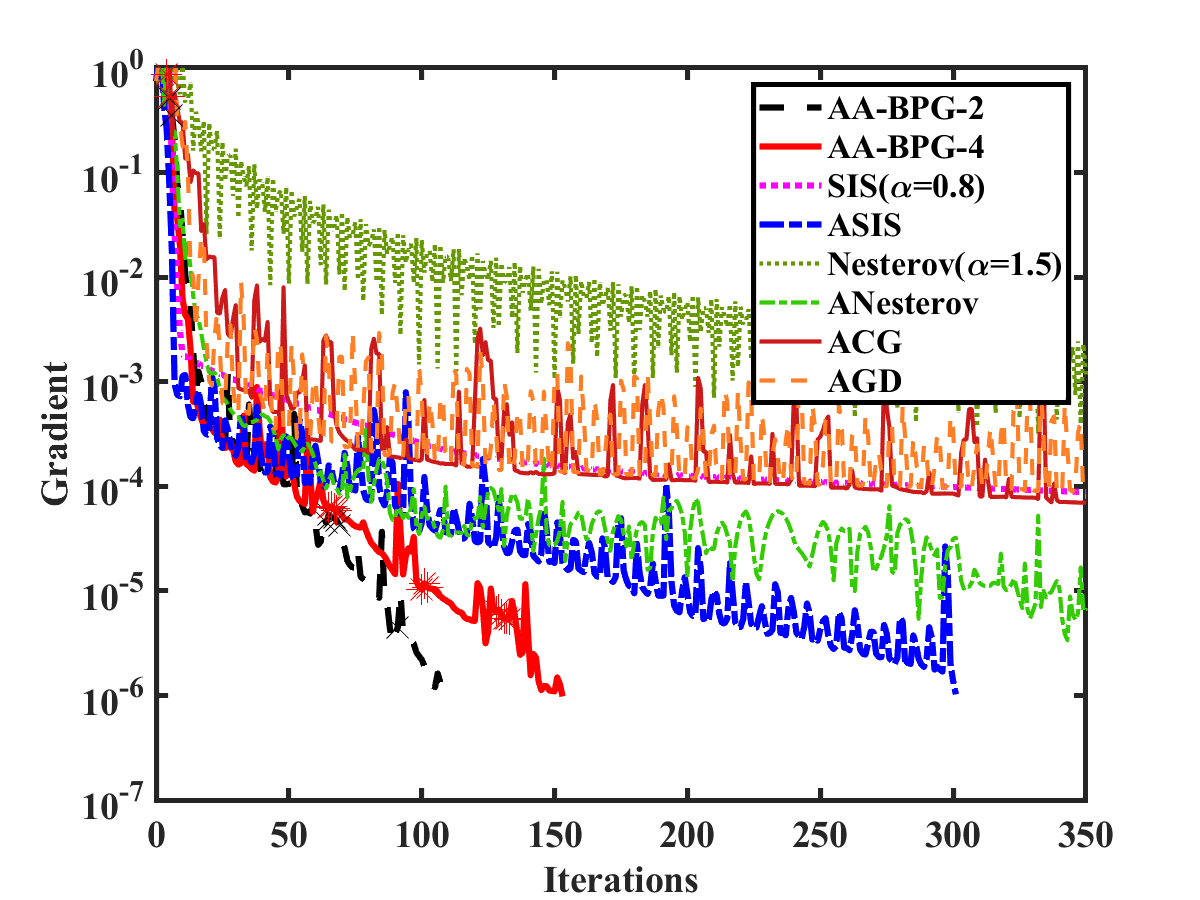}}
	\caption{Comparison of numerical behavior of AA-BPG-2, AA-BPG-4 and the other methods for computing the striped phase on a sphere of radius $R=\sqrt{3660}$.
		The information on these profiles is the same as Fig.\;\ref{fig:BPGk_spot_15}
	}
	\label{fig:BPGk_strip_60}
\end{figure}
It is observed that the AA-BPG-2 and AA-BPG-4 methods cost iterations fewer than $160$ and CPU time less than $50$ seconds, which have much better performance than other existing methods.
The iterative process of these algorithms is given in Fig.\;\ref{fig:BPGk_strip_60}.
Note that the ANesterov algorithm performs poorly here, although it has a faster convergence speed in computing the $60$-spotted phase.
Energy oscillations are observed in the ANesterov and ACG algorithms, while the sequences generated by others keep a monotonic decrease of free energy.
Based on these results, the AA-BPG-2 method is considered the most efficient.
The step sizes of the adaptive methods are illustrated in
Fig.\;\ref{fig:L60_time_step}.
\par
\begin{figure}[!htb]
	\centering
	\includegraphics[width=74mm,height=60mm]{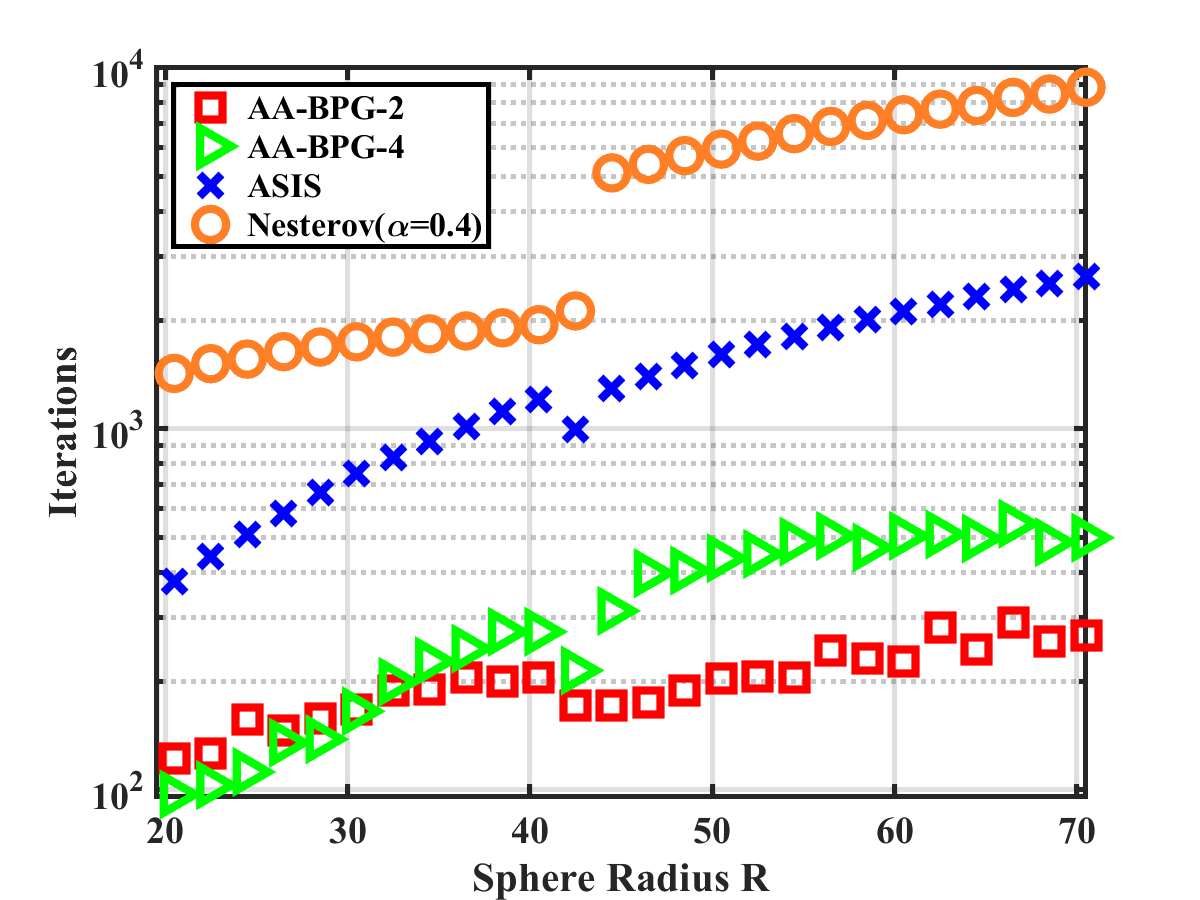}
	\hspace{0.5cm}	
	\includegraphics[width=74mm,height=60mm]{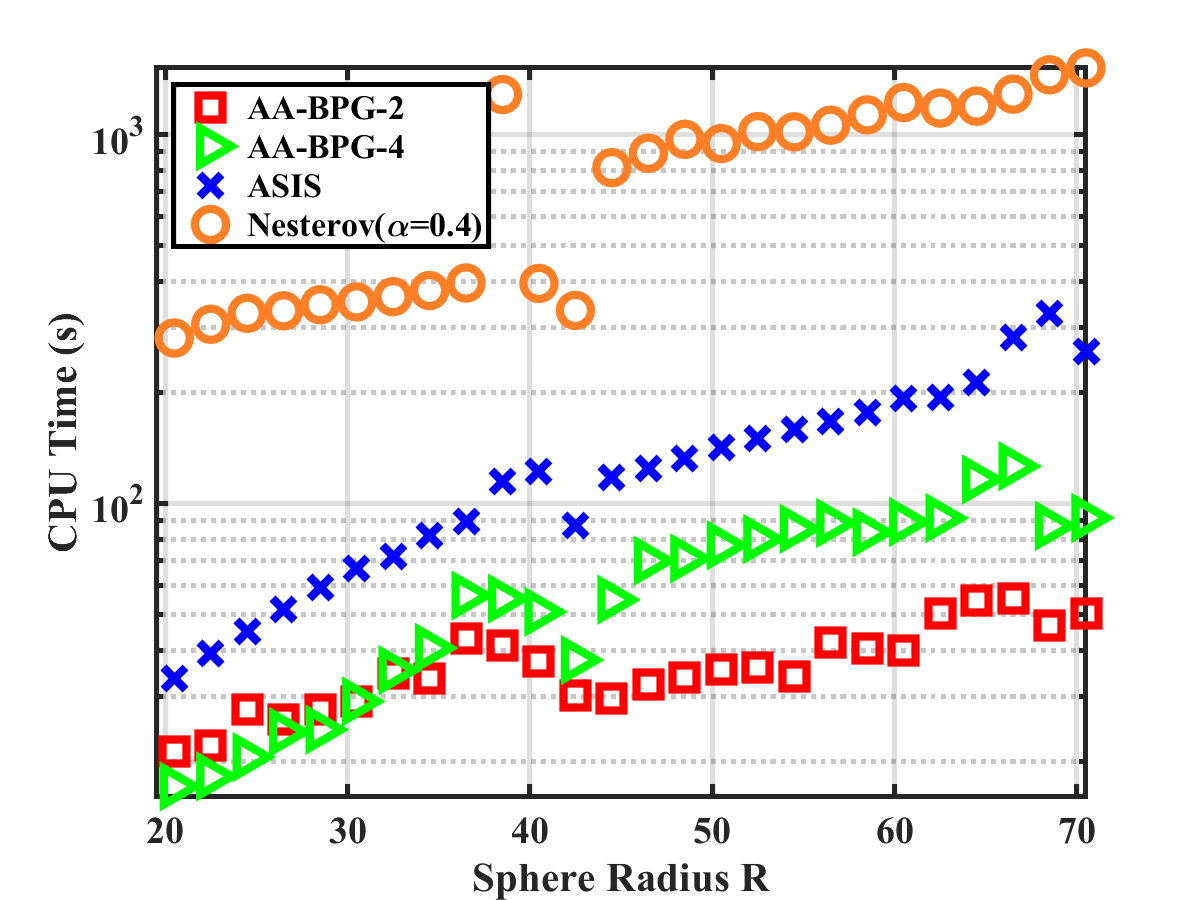}
	\caption{Required iterations and CPU time of the AA-BPG-2, AA-BPG-4, Nesterov and ASIS methods for the stationary striped phase
		with different sphere radius $R$ when $\xi=1.0$, $\epsilon=-0.8$,
		$\lambda=0.0$. Here $\alpha_{0}=0.4$, $\alpha_{\min}=1.0\times 10^{-4}$ and
		$\alpha_{\max}=1.0$ 
	}
	\label{fig:ABPG_striped}
\end{figure}

We further compare the performance of the AA-BPG-2, AA-BPG-4, ASIS and Nesterov methods by computing other striped phases.
We also use the PMA method to choose the initial configurations and sphere radius.
Here we fix $\alpha=0.4$ for the Nesterov method.
Figure\;\ref{fig:ABPG_striped} presents their performance over a sequence of sphere radius $R$.
These profiles again show that the AA-BPG-2 method performs the best. 

\par
\TR{Finally, we take AA-BPG-2/4 as the example to demonstrate the robustness of the proposed optimization algorithms for different striped phases at different model parameters $\xi$, $\epsilon$, $\lambda$.
	Table \ref{tab:AA_BPG_L_60_para} shows the corresponding convergent results.
	Here, we fix $\xi=1.0$, $\alpha_{\min}=0.01$, $\alpha_{\max}=5.0$, but $\alpha_{0}=0.8$ for ASIS and  $\alpha_{0}=0.5$ for AA-BPG methods.
	From the table, we can see that compared to ASIS, AA-BPG methods have better computational efficiency in calculating different striped structures at varying $\epsilon$ and $\lambda$.
}

\begin{table} [!htb] 
	\TR{
		\caption{
			Comparisons of numerical behevior of AA-BPG-2, AA-BPG-4 and ASIS for computing a stiped phase at varying model parameters $\epsilon$, $\lambda$.
		}  
		\begin{center}
			\begin{tabular}{|p{0.8cm}|p{0.8cm}|p{1.2cm}|p{2.0cm}|p{1.8cm}|p{2.2cm}|p{3.8cm}|}  
				\hline  
				$\epsilon$ & $\lambda$ & R	&Method &Iteration & CPU time (s) & Equilibrium energy $E_{s}$   \\  
				\hline  
				\multirow{3}*{-0.9} & \multirow{3}*{0} & \multirow{3}*{$\sqrt{3660}$} & AA-BPG-2 &111 & 23.13 & -2.8647889426\\ 
				~ &~ & ~ & AA-BPG-4 &160 & 39.13 & -2.8647889426\\
				~&~ & ~ & ASIS & 339 & 61.81 & -2.8647889426\\
				\hline
				\hline
				\multirow{3}*{-1} & \multirow{3}*{0.01} & \multirow{3}*{$\sqrt{3660}$} & AA-BPG-2 &160 & 49.40 & -3.5375809017\\
				~ &~ & ~ & AA-BPG-4 &138 & 41.47 & -3.5375809017\\
				~ &~ & ~& ASIS &367 & 82.16 & -3.5375809017\\
				\hline
				\hline
				\multirow{3}*{-0.8} & \multirow{3}*{0.05} & \multirow{3}*{$\sqrt{3660}$} & AA-BPG-2 &184 & 41.46 & -2.2636412022\\
				~ &~ & ~ & AA-BPG-4 &187 & 44.62 & -2.2636412022\\
				~ &~ & ~ & ASIS &345 & 68.19 & -2.2636412022\\
				\hline
				\hline
				\multirow{3}*{-0.9} & \multirow{3}*{0.05} & \multirow{3}*{$\sqrt{6480}$} & AA-BPG-2 &171 & 35.07 & -2.8674146362\\
				~ &~ & ~  & AA-BPG-4 &130 & 29.44 & -2.8674146362\\
				~ &~ & ~  & ASIS &394 &75.21 & -2.8674146362\\
				\hline
				\hline
				\multirow{3}*{-1} & \multirow{3}*{0.05} & \multirow{3}*{$\sqrt{6480}$} & AA-BPG-2 &160 & 33.60 & -3.5404666645\\
				~ &~ & ~  & AA-BPG-4 &174 & 40.26 & -3.5404666645\\
				~ &~ & ~  & ASIS &552 &111 & -3.5404666645\\
				\hline
				\hline
				\multirow{3}*{-0.7} & \multirow{3}*{0.015} & \multirow{3}*{$\sqrt{6480}$} & AA-BPG-2 & 279& 65.42 & -1.7339426725\\
				~ &~ & ~  & AA-BPG-4 &149 & 57.00  & -1.7339426725\\
				~ &~ & ~  & ASIS &346 &141.14 & -1.7339426725\\
				\hline
			\end{tabular}
		\end{center}
		\label{tab:AA_BPG_L_60_para} 
	}
\end{table}

%% file: Conclusion.tex
\section{Conclusion}
\label{sec:Conc}
In this paper, we have developed efficient numerical methods to compute stationary states of the spherical LB model.
Instead of solving the gradient flow equation, we compute stationary states of the free energy function directly using optimization algorithms based on the discrete spherical harmonic expansion.
The developed optimization algorithms include the standard AGD, ACG, Nesterov, ANesterov and AA-BPG methods.
Numerical experiments on different striped and spotted
phases show that the AA-BPG and Nesterov methods significantly reduce the number of iterations and the computational time required for convergence.
Furthermore, we have proposed the PMA method to estimate good initial values to accelerate the convergence to stationary structures.
It indicates that good initial states can be constructed by the relations between symmetric subgroups of $O(3)$ and spherical harmonics of degree $\ell$, and that the sphere radius satisfies $\sqrt{\ell(\ell+1)}$.
Extensive results validate the effectiveness of the proposed approach in accurately and efficiently finding stationary structures with desired symmetry.
In the future, we will extend our methods to explore richer phase behavior of ordered structures on the spherical surface, such as phase transitions.

\section*{Acknowledgments}
G. Ji is partially supported by the National Natural Science Foundation of China (12471363).
K. Jiang is supported in part by the National Key R\&D Program of China (2023YFA1008802), the National Natural Science Foundation of China (12171412), the Natural Science Foundation for Distinguished Young Scholars of Hunan Province (2021JJ10037),
the Science and Technology Innovation Program of Hunan Province (2024RC1052).
We are also grateful to the High Performance Computing Platform of Xiangtan University for partial support of this work.

%% file: sphere_LB.bbl
\begin{thebibliography}{99}
	
	\bibitem{1950_LG}L. D. Landau and V. L. Ginzburg, On the theory of superconductivity, Zh. Eksp. Teor. Fiz., 20 (1950), 1064–1082.
	\bibitem{1975Brazovskii}S. Brazovskii, Phase transition of an isotropic system to a nonuniform state, J. Exp. Theor. Phys., 41 (1975), 85–89.
	\bibitem{1977_SH}J. Swift and P. Hohenberg, Hydrodynamic fluctuations at the convective instability, Phys. Rev. A., 15 (1977), 319–328.
	\bibitem{Zhang2008}P. Zhang and X. Zhang, An efficient numerical method of Landau–Brazovskii model, J. Comput. Phys., 227 (2008), 5859–5870.
	\bibitem{Shi2003}R. A. Wickham, A.-C. Shi and Z. Wang, Nucleation of stable cylinders from a metastable lamellar phase in a diblock copolymer melt, J. Chem. Phys., 118 (2003), 10293–10305.
	\bibitem{Zhanghao2014}H. Zhang, K. Jiang and P. Zhang, Dynamic transitions for Landau-Brazovskii model, Discrete. Cont. Dyn-B., 19 (2014), 607–627.
	\bibitem{pattern_forming_2003}P. C. Matthews, Pattern formation on a sphere. Phys. Rev. E., 67 (2003), 036206.
	\bibitem{spiral_symmetry_2010}R. Sigrist and P. C. Matthews, Symmetric spiral patterns on spheres, SIAM J. Appl. Dyn. Syst., 10 (2011), 1177–1211.
	\bibitem{Zhang_2014}L. Zhang, L. Wang and J. Lin, Defect structures and ordering behaviours of diblock copolymers self-assembling on spherical substrates, Soft Matter., 10 (2014), 6713–6721.
	\bibitem{SLB_Sanjay_2017}S. Dharmavaram, F. Xie, W. Klug, J. Rudnick and R. Bruinsma, Orientational phase transitions and the assembly of viral capsids, Phys. Rev. E., 95 (2017), 062402.
	\bibitem{2002_Elder_SH}K. R. Elder, M. Katakowski, M. Haataja and M. Grant, Modeling elasticity in crystal growth, Phys. Rev. Lett., 88 (2002), 245701.
	\bibitem{2010_Backofen_Particles}R. Backofen, A. Voigt, Axel and T. Witkowski, Particles on curved surfaces: A dynamic approach by a phase-field-crystal model, Phys. Rev. E., 81 (2010), 025701.
	\bibitem{2011_Backofen_Continuous}R. Backofen, M. Gr\"{a}f, D. Potts, S. Praetorius, A. Voigt and T. Witkowski, A continuous approach to discrete ordering on $\mathbb{S}^2$, Multiscale. Model. Sim., 9 (2011), 314-334.
	\bibitem{2011_Aland_Continuum}S. Aland, J. Lowengrub and A. Voigt, A continuum model of colloid-stabilized interfaces, Phys. Fluids., 23 (2011), 062103.
	\bibitem{2012_Aland_Particles}S. Aland, J. Lowengrub and A. Voigt, Particles at fluid-fluid interfaces: A new Navier-Stokes-Cahn-Hilliard surface- phase-field-crystal model, Phys. Rev. E. 86 (2012), 046321.
	\bibitem{2012_Aland_Bucking}S. Aland, A. R\"{a}tz, M. R\"{o}ger and A. Voigt, Buckling instability of viral capsids - A continuum approach, Multiscale. Model. Sim., 10 (2012), 82-110.
	Multiscale. Model. Sim., 10 (2012), 82-110.
	\bibitem{2018_Praetorius_Active}S. Praetorius, A. Voigt, R.Wittkowski and H. L\"owen, Active crystals on a sphere, Phys. Rev. E., 97 (2018), 052615.
	\bibitem{1987Fredrickson}G. H. Fredrickson and E. Helfand, Fluctuation effects in the theory of microphase separation in block copolymers, J. Chem. Phys., 87 (1987), 697–705.
	\bibitem{SIS_Shen_2010}J. Shen  and X. Yang, Numerical approximations of Allen-Cahn and Cahn-Hilliard equations, Discrete. Cont. Dyn-A., 28 (2010), 1669–1691.
	\bibitem{CN_Feng_2013}X. Feng, T. Tang and J. Yang, Stabilized Crank-Nicolson/Adams-Bashforth schemes for phase field models, E. ASIAM. J. Appl. Math., 3 (2013), 59–80.
	\bibitem{CSplitting_Vignal_2015}P. Vignal, L. Dalcin, D. L. Brown, N. Collier and V. M. Calo, An energy-stable convex splitting for the phase-field crystal equation, Comput. Struct., 158 (2015), 355–368.
	\bibitem{CSplitting_Lee_2016}J. Shin, H. G. Lee and J.-Y. Lee, First and second order numerical methods based on a new convex splitting for phase-field crystal equation, J. Comput. Phys., 327 (2016), 519–542.
	\bibitem{OSplitting_Lee_2015}H. G. Lee, J. Shin and J.-Y. Lee, First and second order operator splitting methods for the phase field crystal equation, J. Comput. Phys., 299 (2015), 82–91.
	\bibitem{OSplitting_Zhai_2021}S. Zhai, Z. Weng, X. Feng and Y. He, Stability and error estimate of the operator splitting method for the phase field crystal equation, J. Sci. Comput., 86 (2021), 1–23.
	\bibitem{IEQ_Yang_2016}X. Yang, Linear first and second-order unconditionally energy stable numerical schemes for the phase field model of homopolymer blends, J. Comput. Phys., 327 (2016), 294–316.
	\bibitem{CN_Yang_2017}X. Yang and D. Han, Linearly first- and second-order unconditionally energy stable schemes for the phase field crystal model, J. Comput. Phys., 330 (2017), 1116–1134.
	\bibitem{SAV_Shen_2019}J. Shen, J. Xu and J. Yang, A new class of efficient and robust energy stable schemes for gradient flows, SIAM Rev., 61 (2019), 474–506.
	\bibitem{2005_Tang_Phase}P. Tang, F. Qiu, H. Zhang and Y. Yang,
	Phase separation patterns for diblock copolymers on spherical surfaces: A finite volume method, Phys. Rev. E., 72 (2005), 016710.
	\bibitem{2022_Yang_FDM}J. Yang, J. Wang and Z. Tan, A simple and practical finite difference method for the phase-field crystal model with a strong nonlinear vacancy potential on 3D surfaces, Comput. Math. Appl., 121 (2022), 131-144.
	\bibitem{Luo2018SPH_mini}Y. Luo and L. Maibaum, Phase diagrams of multicomponent lipid vesicles: Effects of finite size and spherical geometry, J. Chem. Phys., 149 (2018), 174901.
	\bibitem{JiangNBPG}K. Jiang, W. Si, C. Chen and C. Bao, Efficient numerical methods for computing the stationary states of phase field crystal models, SIAM J. Sci. Comput., 42 (2020), 1350–1377.
	\bibitem{2021JiangBPG}C. Bao, C. Chen and K. Jiang, An adaptive block bregman proximal gradient method for computing stationary states of multicomponent phase-field crystal model, CSIAM Trans. Appl. Math. 3 (2022), 133–171.
	\bibitem{2012Spherical}K. Atkinson and W. Han, Spherical harmonics and approximations on the unit sphere: An introduction, Springer Berlin Heidelberg, 2012.
	\bibitem{Schulten1975wigner}K. Schulten and R. G. Gordon, Exact recursive evaluation of 3j- and 6j-coefficients for quantum-mechanical coupling of angular momenta,
	J. Math. Phys., 16 (1975), 1961–1970.
	\bibitem{Dahlen_wigner_1998}F. Dahlen and J. Tromp, Theoretical global seismology, Princeton University Press, 1988.
	\bibitem{shtns_Ishioka_2018}K. Ishioka, A new recurrence formula for efficient computation of spherical harmonic transform, J. Me. Teorol. Soc. JPN., 96 (2018), 241–249.
	\bibitem{BB_Barzilai_1988}J. Barzilai and J. M. Borwein, Two-point step size gradient methods, IMA. J. Numer. Anal., 8 (1988), 141–148.
	\bibitem{Stability_LP_Jiang}K. Jiang, J. Tong, P. Zhang and A.-C. Shi, Stability of two-dimensional soft quasicrystals in systems with two length scales, Phys. Rev.
	E., 92 (2015), 042159.
	\bibitem{2018_Sanjay_symmetry}S. Dharmavaram and T. Healey, Direct construction of symmetry-breaking directions in bifurcation problems with spherical symmetry,
	Discrete. Cont. Dyn-S., 12 (2018), 1669–1684.
	
\end{thebibliography}
